\documentclass[12pt]{article}

\usepackage{a4wide,amssymb,amsmath,latexsym,amsthm}
\usepackage[all]{xy}

\newtheorem{theorem}{Theorem}[section]
\newtheorem{definition}[theorem]{Definition}
\newtheorem{lemma}[theorem]{Lemma}

\newtheorem{claim}[theorem]{Claim}

\newtheorem{proposition}[theorem]{Proposition}
\newtheorem{corollary}[theorem]{Corollary}

\newenvironment{prf}[1]{\trivlist 
\item[\hskip \labelsep{\it 
#1.\hspace*{.3em}}]}{~\hspace{\fill}~$\square$\endtrivlist} 

\newcommand{\spec}[1]{ \mathrm{Spec} (#1)}

\newcommand{\xq}{\mathbb{Q}}
\newcommand{\QQ}{\mathbb{Q}}

\newcommand{\ZZ}{\mathbb{Z}}
\newcommand{\xz}{\mathbb{Z}}
\newcommand{\xf}{\mathbb{F}}
\newcommand{\FF}{\mathbb{F}}
\newcommand{\xg}{\mathbb{G}}
\newcommand{\xp}{\mathbb{P}}

\newcommand{\pol}{\mathcal{L}}
\newcommand{\polp}{\mathcal{L}^{(p)}}
\newcommand{\emm}{\mathcal{M}}
\newcommand{\ab}[3]{#1:#2\rightarrow#3}
\newcommand{\hsp}{\hspace{5pt}}

\begin{document}

\title{Galois theory of the canonical theta structure}
\author{Robert Carls}

\maketitle

\begin{abstract}
\noindent
In this article we give a Galois-theoretic characterization
of the canonical theta structure.
The Galois property of the canonical theta structure
translates into certain $p$-adic theta relations
which are satisfied by the canonical theta null point
of the canonical lift.
As an application we give a purely algebraic proof of some $2$-adic
theta identities which describe the set of theta null points
of the canonical lifts of ordinary abelian varieties in
characteristic $2$.
The latter theta relations are suitable for explicit canonical lifting.
\end{abstract}

\noindent
{\small
\begin{center}
\begin{tabular}{ll}
KEYWORDS: & Canonical Lifting, $p$-Adic Theta Functions \\
AMS MSC: & 11G15, 14K25
\end{tabular}
\end{center}
}

\section{Introduction}

The main result of this article is a proof of
the characterizing Galois property
of the canonical theta structure.
The results of this article complement the discussion in
\cite[$\S$2.2]{ckl08}.
The proofs given there are valid for theta
structures which are `away from $p$'. In this article we deal with
the inseparable case, i.e we focus on theta structures
centered `at $p$'.
Our main result forms an important ingredient of the method
for the computation of canonical lifts that is
given in \cite{cl09}.
\newline\indent
In this section we briefly review our main theorem
and give an application to the characterization
of the theta null points of $2$-adic
canonical lifts.
In the remaining sections of the introduction
we show how the canonical theta structure
can be used to give a theoretical
foundation to Mestre's AGM point counting
algorithm.
\newline\indent
Assume that we are given the canonical lift $A$
of an ordinary abelian
variety over a perfect field of characteristic $p>0$
and that $A$ admits an ample symmetric line bundle $\pol$
of degree $1$.
It is proved in \cite{ca07} that there exists a canonical theta 
structure $\Theta$ for the $p$-th power line bundle $\pol^p$.
We note that the canonical theta structure depends on the choice of
a Lagrangian structure on the special fiber.
The theta null point with respect to the canonical theta structure
gives an arithmetic invariant which is attached to the triple $(A,\pol,\Theta)$.
\newline\indent 
It is classically known that there exists a canonical lift
$F:A \rightarrow A^{(p)}$
of the relative $p$-Frobenius morphism.
Assume that there exists a $p$-adic lift $\sigma$
of the $p$-th power Frobenius automorphism of the residue field.
In the following we describe the interplay of the action of $F$ and $\sigma$
on the arithmetic invariant of the canonical lift which is
induced by the canonical theta structure.
First, consider the following construction using the Frobenius lift $F$.
One can prove that there exists an induced ample symmetric
line bundle $\pol^{(p)}$ of degree $1$
on $A^{(p)}$ such that $F^* \pol^{(p)} \cong \pol^p$.
The $p$-th power of the line bundle $\pol^{(p)}$ carries a natural theta structure $\Theta^{(p)}$
which one obtains by canonically descending the canonical theta structure $\Theta$
along the Frobenius lift $F$.
As a consequence, one gets an arithmetic invariant of $A^{(p)}$ by taking the
theta null point with respect to the triple $\big( A^{(p)},\big(
\pol^{(p)} \big)^p,\Theta^{(p)} \big)$.
Another way of defining an arithmetic invariant of $A^{(p)}$
is by twisting the triple $(A,\pol,\Theta)$ along the automorphism $\sigma$.
We denote the resulting triple by $\big(
A^{(\sigma)},\big( \pol^{(\sigma)} \big)^p,\Theta^{(\sigma)} \big)$.
Since $A$ is assumed to be a canonical lift, it follows that
$A^{(p)}$ and $A^{(\sigma)}$, and the $p$-th powers of the line bundles
$\pol^{(p)}$ and $\pol^{(\sigma)}$, are isomorphic.
In this article we prove the following equality.
\begin{theorem}
\label{canthetapb}
Up to canonical isomorphism one has
\begin{eqnarray}
\label{equa}
\Theta^{(p)} = \Theta^{(\sigma)}.
\end{eqnarray}
\end{theorem}
\noindent
We have stated the above theorem in a slightly informal way to avoid
technicalities.
For more details regarding the formalism
and a proof see Theorem \ref{cantwisttheta}
of Section \ref{galoisprop}.
The above theorem tells us that the two constructions from above, namely the descent along $F$ and the twisting by $\sigma$, define the
same arithmetic invariant of $A^{(p)}=A^{(\sigma)}$.
The equality (\ref{equa})
implicitly gives a set of special $p$-adic theta relations
which hold exclusively for the canonical theta null point, i.e.
the theta null point with respect to the canonical theta structure.
\newline\indent
Using the
Galois-theoretic characterization of the canonical theta
structure in the case of residue characteristic $2$
one can prove some $2$-adic theta relations which describe the set of
theta null points of canonical lifts.
Now let $p=2$ and $d \geq 1$ an integer. Assume that the abelian scheme
$A$ from above is defined over
the unramified extension $\ZZ_q$ of degree $d$ of the $2$-adic integers.
We denote the theta null point induced by the canonical theta structure
$\Theta$ by $( a_u )_{u \in (\ZZ / 2 \ZZ)^g}$.
Combining Mumford's $2$-multiplication formula \cite[$\S$3]{mu66}
with the theta identities coming from the Galois theoretic
characterization
of the canonical theta structure (compare Theorem \ref{canthetapb})
we get the following $2$-adic theta
identities.
\begin{theorem}
\label{2adictheta}
There exists an $\omega \in \ZZ_q^*$ such that for all $u \in (\xz / 2 \xz)^g$
the following equality holds
\begin{eqnarray}
\label{specialcase}
a_u^2 = \omega 
\sum_{v \in (\xz / 2 \xz)^g} a^{\sigma}_{ u + v }  a^{\sigma}_v.
\end{eqnarray}
\end{theorem}
\noindent
This statement is a trivial corollary of Theorem \ref{lst} which
is proved in Section \ref{main}.
We illustrate the above results by examples
in Section \ref{ccft}.
We remark that the theta relations (\ref{specialcase})
form an essential ingredient of
the quasi-quadratic point counting algorithm
that is given in \cite{ll06}.
By applying a multivariate Hensel lifting algorithm
to the theta identities (\ref{specialcase}) one obtains
an explicit method for higher dimensional canonical lifting.
Practical aspects of the lifting of the canonical theta null point
on the basis of the equations (\ref{specialcase}) are
discussed in \cite[$\S$4.1]{saga}.
We would like to point out that the algorithm given there
produces completely integral recursions locally at $2$, which contrasts 
with the classical analytic AGM recursions used by Mestre (see Section
\ref{workofmestre}), that require 
division by $2^g$ in the $2$-adic Witt ring. This in itself yields a more 
efficient and stable algorithm.

\subsection{On Mestre's generalized $p$-adic AGM sequence}
\label{onagm}

The results of the preceding section can be used to
give a theoretical foundation to Mestre's
AGM algorithm.
In the year 2002 Mestre proposed an algorithm for point counting
on ordinary hyperelliptic curves over finite fields
of characteristic $2$ which is based on the computation
of the generalized $2$-adic \emph{arithmetic geometric mean} (AGM)
sequence (see the informal notes \cite{me02}).
The latter algorithm extends earlier results of his on
the algorithmic aspects of the $p$-adic Gaussian AGM (see
\cite{hm89} \cite{me00}).
Lercier and Lubicz have shown in \cite{ll06} how
Mestre's algorithm can be modified such that its asymptotic
complexity becomes quasi-quadratic in the degree of the finite field.

\subsubsection{Elementary properties of the AGM sequence}
\label{workofmestre}

In this section we recall
a result of Mestre about the convergence of the
generalized $p$-adic AGM
sequence.
By lack of a suitable reference we also provide a proof.
Consider the multi-valued sequence of vectors
$\big( x^{(n)}_u \big)_{ u \in (\mathbb{Z} / 2 \mathbb{Z})^g}$
which is given by
the recursive law
\begin{eqnarray}
\label{gagm}
x^{(n+1)}_u = \frac{1}{2^g} \sum_{v \in (\mathbb{Z} / 2 \mathbb{Z})^g} \sqrt{x^{(n)}_{u+v} x^{(n)}_v},
\end{eqnarray}
where $u \in (\mathbb{Z} / 2 \mathbb{Z})^g$.
For $g=1$ the above formulas resemble
the classical AGM formulas of Gauss \cite{co84} \cite{bm88}.
The case $g=2$ was first studied by Borchardt \cite{bo1879}.
Under suitable assumptions, that are made precise in the following
theorem, the sequence (\ref{gagm})
is well defined over the $p$-adic numbers.
\newline\indent
Let $g \geq 1$ and let $\FF_q$ be a finite field of
characteristic $p>0$ which has
$q$ elements.
We denote the Witt vectors with values in $\FF_q$
by $\ZZ_q$.
\begin{lemma}[Mestre]
If one starts with
$x_u^{(0)} \in 1 + 2^{g+1}p \ZZ_q$
and if one chooses all square roots
$\equiv 1 \bmod 2p$ then the recursion {\rm (\ref{gagm})} gives a well-defined sequence
in the ring $\ZZ_q$.
Under the assumption
\begin{eqnarray}
\label{convergence}
x_u^{(0)} \in 1 + 2^{g+2} p \ZZ_q, \quad \forall u \in (\ZZ / 2 \ZZ)^g
\end{eqnarray}
the generalized AGM sequence {\rm (\ref{gagm})}
converges coordinatewise, and one has
\begin{eqnarray}
\label{seq}
\lim_{ n \rightarrow \infty} x_u^{(n)}=\lim_{ n \rightarrow \infty}
x_0^{(n)}, \quad \forall u \in (\ZZ / 2 \ZZ)^g.
\end{eqnarray}
\end{lemma}
\begin{proof}
It is a basic fact that the squares in the ring $\xz_q$ that are congruent $1$
modulo $2p$ are given by the
subset $1+4p \xz_q$. 
In the following let $v_p:\xz_q^* \rightarrow \xz$ denote an exponential
valuation such that $v_p(p)=1$.
We set $|x|_p=p^{-v_p(x)}$ for $x \in \xz_q^*$.
If $x \in 1+4p \xz_q$, then there exists a unique
square root $\sqrt{x}$ of $x$ which is congruent $1$ modulo $2p$.
Moreover, if $p=2$
then $v_2 \big( \sqrt{x}-1 \big)=v_2(x-1)-1$.
Otherwise, if $p>2$ then $v_p(x-1)=v_p \big( \sqrt{x}-1 \big)$.
\newline\indent
Now assume that $x_u^{(0)} \in 1 + 2^{g+1}p \ZZ_q$ for
all $u \in (\xz / 2 \xz)^g$.
By the above discussion the value $x_u^{(0)}$ is a square.
Let $y_u^{(0)}$ be the unique square root of $x_u^{(0)}$ such that
$y_u^{(0)} \equiv 1 \bmod 2p$.
Let now $x_u^{(1)}$ be defined by the formulas (\ref{gagm}), where the
square roots are chosen $\equiv 1 \bmod 2p$.
The formula (\ref{gagm}) reads as
\begin{eqnarray}
\label{modagm}
x_u^{(1)} =\frac{1}{2^g}
\sum_{v \in (\xz / 2 \xz)^g} y_{u+v}^{(0)} \cdot y_v^{(0)},
\quad u \in (\ZZ / 2 \ZZ)^g.
\end{eqnarray}
Let $w \in (\xz / 2 \xz)^g$ and let $S_w$ be a subset of
$(\xz / 2 \xz)^g$ which is maximal with respect to the
property that for all $v \in S_w$ one has $v+w \not\in S_w$.
By assumption we get
\begin{eqnarray}
\label{tail}
\lefteqn{ \left| x_0^{(1)} - x_w^{(1)}
\right|_p
\stackrel{(\ref{modagm})}{=} \frac{1}{|2|_p^g} \cdot \left| \sum_{v \in (\xz / 2 \xz)^g}
\left[ \Big( y_v^{(0)} \Big)^2-y_{w+v}^{(0)} \cdot y_{v}^{(0)}
\right] \right|_p} \\
\nonumber && = \frac{1}{|2|_p^g} \cdot \left| \sum_{ v \in S_w}
\Big( y_v^{(0)}-y_{w+v}^{(0)} \Big)^2 \right|_p
\leq \frac{1}{|2|_p^g} \cdot \mathrm{max}_{ v \in (\xz / 2 \xz)^g}
\left| y_v^{(0)} - y_{w+v}^{(0)} \right|_p^2
\end{eqnarray}
It follows that
\begin{eqnarray*}
\mathrm{max}_{ v \in (\xz / 2 \xz)^g} \big| x_0^{(1)}-x_v^{(1)} \big|_p
\leq \frac{1}{|2|_p^g} \cdot \mathrm{max}_{ v \in (\xz / 2 \xz)^g}
\left| y_0^{(0)} - y_{v}^{(0)} \right|^2_p
\end{eqnarray*}
The latter equation implies that $x_u^{(1)} \in 1+2^{g+1}p \xz_q$.
This shows that the AGM sequence is well-defined.
\newline\indent
Now let $y_u^{(1)}$ be the unique square root of $x_u^{(1)}$ such that
$y_u^{(1)} \equiv 1 \bmod 2p$.
One computes
\begin{eqnarray}
\label{head}
|2|_p \cdot \big| y_0^{(1)}-y_w^{(1)} \big|_p
= \big| y_0^{(1)}+y_w^{(1)} \big|_p \cdot \big| y_0^{(1)}-y_w^{(1)}
\big|_p=\big| x_0^{(1)} - x_w^{(1)}
\big|_p
\end{eqnarray}
Suppose that we have $x_u^{(0)} \in 1+2^{g+2}p \xz_q$.
Under this assumption, by combining (\ref{head}) and (\ref{tail}), we get
\begin{eqnarray*}
\mathrm{max}_{ v \in (\xz / 2 \xz)^g} \big| y_0^{(1)}-y_v^{(1)} \big|_p
\leq \frac{1}{2} \cdot \mathrm{max}_{ v \in (\xz / 2 \xz)^g}
\left| y_0^{(0)} - y_{v}^{(0)} \right|_p
\end{eqnarray*}
We conclude by induction that
the sequence of vectors $\big( y_0^{(n)} - y_u^{(n)} \big)$,
and hence also the sequence $\big( x_0^{(n)} - x_u^{(n)} \big)$, converge
to zero for all $u \in (\ZZ / 2 \ZZ)^g$.
The coordinatewise convergence now follows from the estimate
\begin{eqnarray*}
\lefteqn{\left| x_0^{(n+1)}-x_0^{(n)} \right|_p
= \frac{1}{|2|_p^g} \cdot
\left| \left( \sum_{v \in (\ZZ / 2 \ZZ)^g} x_v^{(n)} \right)
-2^g \cdot x_0^{(n)} \right|_p} \\
&& \leq \frac{1}{|2|_p^g} \cdot \mathrm{max}_{v \in (\xz / 2 \xz)^g}
 \left| x_0^{(n)} -x_v^{(n)} \right|_p
\end{eqnarray*}
This finishes the proof of the lemma.
\end{proof}

\subsection{$2$-Adic approximation of canonical lifts}
\label{agmscheme}

In this section we give a scheme theoretic counterpart
of the generalized AGM sequence.
We use the notation of the previous section.
Let $q=2^d$ where $d \geq 1$ is an integer.
As before, let $\xz_q$ denote the Witt vectors
with values in a finite field $\xf_q$ with $q$ elements.
Suppose that
$x_u^{(0)} \in 1 + 2^{g+2} \ZZ_q$,
but $x_u^{(0)} \not \in 1 + 2^{g+3} \ZZ_q$.
In this case the generalized AGM sequence (\ref{gagm}) is well-defined
but does not converge coordinatewise.
Mestre observed that,
under suitable initial conditions,
the sequence of $\xz_q$-valued vectors
\begin{eqnarray}
\label{fractionalagm}
\left( \frac{x_u^{(dn)}}{x_0^{(dn)}} \right)_{u \in (\ZZ / 2 \ZZ)^g}
\end{eqnarray}
converges for $n \rightarrow \infty$.
In other words, the sequence of
projective points $\big( x_u^{(dn)} \big)$ converges
to a unique projective point of the projective space $\xp^{2^g-1}$.
This convergence forms
a crucial ingredient of Mestre's AGM algorithm.
In the following we will give a theoretical reason
for this phenomenon.
\newline\indent
We define a scheme theoretic counterpart of the AGM sequence as follows.
Let $A$ be an abelian scheme of relative dimension $g$
over $\xz_q$ which has ordinary reduction. It is well-known
that there exists a unique isogeny $F:A \rightarrow A^{(2)}$
which lifts the relative Frobenius on the special fiber.
One has $\mathrm{Ker}(F)=A[2]^{\mathrm{loc}}$, where the latter
denotes the connected component of the $2$-torsion subgroup
$A[2]$.
Successively dividing out the connected component of the
$2$-torsion gives a sequence of
abelian schemes and Frobenius lifts
\begin{eqnarray}
\label{stagm}
A \stackrel{F}{\rightarrow} A^{(2)}
\stackrel{F}{\rightarrow} A^{(4)} \stackrel{F}{\rightarrow} A^{(8)} \rightarrow \ldots
\end{eqnarray}
Let us fix the following notation.
If $n \geq 1$ and $B$ is an abelian scheme over $\xz_q$, then we denote its
reduction modulo $2^n$ by $B_{2^n}$.
Furthermore, we set $A^{(1)}=A$.
Obviously, there exists a canonical isomorphism
$A_2 \stackrel{\sim}{\rightarrow} A^{(2^{dn})}_2$
for all $n \geq 0$.
Let $\tilde{A}^{(2^n)}$ denote the canonical lift of $A^{(2^n)}_2$.
\begin{theorem}
\label{stconvergence}
The sequence $A^{(q^n)}$ converges $2$-adically
to the abelian scheme $\tilde{A}^{(1)}$.
More precisely, for all $n \geq 0$ there exists an
isomorphism $A^{(2^n)}_{2^{n+1}}
\stackrel{\sim}{\rightarrow} \tilde{A}^{(2^n)}_{2^{n+1}}$.
\end{theorem}
\noindent
A more general version
of the above theorem with detailed proof
can be found in \cite[Th.2.1.2]{ca04}.
\newline\indent
Now assume that $A$ comes equipped with an ample symmetric
line bundle $\pol$ of degree $1$.
Then by \cite[Th.5.1]{ca07} for all $n \geq 1$ the line bundle $\pol^{2^n}$
canonically descends along the isogeny $F^n$
to an ample symmetric line bundle $\pol^{(2^n)}$
of degree $1$ on $A^{(2^n)}$.
We set $\pol^{(1)}=\pol$.
Assume that we are given an isomorphism
\begin{eqnarray}
\label{pptriv}
(\xz / 2 \xz)^g
\stackrel{\sim}{\rightarrow} A[2]^{\mathrm{et}}
\end{eqnarray}
where $A[2]^{\mathrm{et}}$ denotes the maximal \'etale
quotient of $A[2]$.
By composition with $F^n$ the isomorphism (\ref{pptriv}) induces 
a trivialization
$(\xz / 2 \xz)^g
\stackrel{\sim}{\rightarrow} A^{(2^n)}[2]^{\mathrm{et}}$
for all $n \geq 1$.
By \cite[Th.2.1]{ca07} for every natural number $n \geq 2$
there exists a canonical symmetric theta structure
$\Theta^{(2^n)}$ of type $(\xz / 2 \xz)^g$ for the $2$-nd power of
the line bundle $\pol^{(2^n)}$ which only depends on the
choice of the trivialization (\ref{pptriv}).
To the sequence of triples $ \big(
A^{(2^n)}, \big(\pol^{(2^n)} \big)^2,\Theta^{(2^n)} \big)$
one can associate
a sequence of theta null points $\big( a_u^{(n)} \big)_{u
\in (\xz / 2 \xz)^g}$ in the
projective space $\xp^{2^g-1}_{\xz_q}$.
The line bundle $\pol^{(1)}$ on $A$ induces an ample symmetric line
bundle $\pol^{(1)}_2$ of degree $1$ on the reduction of $A$.
By general theory the line bundle $\pol^{(1)}_2$ lifts
to an ample symmetric
line bundle $\tilde{\pol}^{(1)}$ of degree $1$ on $\tilde{A}^{(1)}$.
We note that the class of the line bundle $\big( \tilde{\pol}^{(1)}
\big)^2$ does
not depend on the chosen lift $\tilde{\pol}^{(1)}$ of $\pol^{(1)}_2$.
By \cite[Th2.2]{ca07} there exists a canonical
theta structure $\tilde{\Theta}^{(1)}$ of type $(\xz / 2 \xz)^g$ for the
$2$-nd power of the
line bundle $\tilde{\pol}^{(1)}$ which only depends on the reduction
of the trivialization (\ref{pptriv}).
We denote the theta null point with respect to the triple
$\big( \tilde{A}^{(1)},\big( \tilde{\pol}^{(1)} \big)^2,\tilde{\Theta}^{(1)} )$
by $\big( \tilde{a}^{(0)}_u \big)$.
\begin{theorem}
\label{breuk}
The sequence of points $\big( a_u^{(dn)} \big)$
converges to the point $\big( \tilde{a}^{(0)}_u \big)$
in the projective space $\xp_{\xz_q}^{2^g-1}$ in the
following sense:
We have $a_0^{(dn)},\tilde{a}^{(0)}_0 \in \xz_q^*$ and
\[
\lim_{n \rightarrow \infty}
\frac{a_u^{(dn)}}{a_0^{(dn)}} = \frac{\tilde{a}^{(0)}_u}{
\tilde{a}^{(0)}_0}
\]
for all $u \in (\xz / 2 \xz)^g$.
\end{theorem}
\begin{proof}
By Theorem \ref{stconvergence} we have an isomorphism
$A_{2^{dn+1}}^{(2^{dn})} \cong \tilde{A}^{(1)}_{2^{dn+1}}$
for all $n \geq 0$.
We claim that $\big( \tilde{\pol}^{(1)}_{2^{dn+1}} \big)^2
\cong \big( \pol^{(2^{dn})}_{2^{dn+1}} \big)^2$.
By writing $\tilde{\pol}^{(1)}_{2^{dn+1}}$
and $\pol^{(2^{dn})}_{2^{dn+1}}$ we indicate that we are working
over the quotient ring of $\xz_q$ modulo $2^{dn+1}$.
By abuse of notation we consider the line bundle
$\pol^{(2^{dn})}_{2^{dn+1}}$ to be defined on the abelian
scheme $\tilde{A}^{(1)}_{2^{dn+1}}$.
To prove our claim it suffices to show that the line bundles
$\tilde{\pol}^{(1)}_{2^{dn+1}}$ and
$\pol^{(2^{dn})}_{2^{dn+1}}$ differ by an element
of $\mathrm{Pic}^0_{\tilde{A}^{(1)}_{2^{dn+1}}}[2]$.
But this follows immediately from the fact that both line bundles
are symmetric and coincide on the special fiber.
Hence our claim is proved.
\newline\indent
By abuse of notation we consider $\Theta^{(2^{dn})}_{2^{dn+1}}$,
the reduction of $\Theta^{(2^{dn})}$ modulo $2^{2d+1}$,
to be defined
on the abelian scheme $\tilde{A}^{(1)}_{2^{dn+1}}$.
We claim that $\Theta^{(2^{dn})}_{2^{dn+1}}
=\tilde{\Theta}^{(1)}_{2^{dn+1}}$.
Our claim follows from the following considerations.
As canonical theta structures, both theta structures are uniquely determined by
the a trivialization, which in both cases equals the isomorphism (\ref{pptriv}).
As a consequence, the equality
of theta structures holds by construction
(for more details see \cite[Th.2.1-2.2]{ca07}). 
\newline\indent
It remains to show that $a_0^{(dn)}$ and $\tilde{a}^{(0)}_0$ are units
in $\xz_q$. By the above discussion we have
$a_0^{(dn)} \equiv \tilde{a}^{(0)}_0 \bmod 2$.
Applying Theorem \ref{2adictheta} we deduce that for all
$0 \not= u \in (\xz / 2 \xz)^g$ we have $\tilde{a}^{(0)}_u
\equiv 0 \bmod 2$. By ampleness we have
$\tilde{a}^{(0)}_0 \not\equiv 0 \bmod 2$.
This completes the proof of the theorem.
\end{proof}
\begin{theorem}
For all $n \geq 2$ there exists a
$\lambda_n \in \xz_q^*$ such that
\begin{eqnarray}
\label{ctrelation}
\Big( a_u^{(n)} \Big)^2= \lambda_n \cdot \sum_{v \in (\xz / 2 \xz)^g}
a_{u+v}^{(n+1)} \cdot a_v^{(n+1)}, \quad u \in (\xz / 2 \xz)^g.
\end{eqnarray}
\end{theorem}
\begin{proof}
The proof is very similar to the one of Theorem \ref{lst},
so we don't give it here.
The most important ingredient of the proof,
namely the $F$-compatibility of the theta structures
$\Theta^{(2^n)}$ and $\Theta^{(2^{n+1})}$,
is implied by slightly modified versions of
Lemma \ref{comp1} and of Lemma \ref{comp2}.
\end{proof}
\noindent
As an obvious consequence of Theorem \ref{2adictheta} the points
$\big( a_u^{(n)} \big)$ lie in the neighborhood
of the point $(1,0,\ldots,0)$ while the points
in Mestre's AGM sequence are contained in a $2$-adic disc
around the point $(1, \ldots,1)$.
In the following we will give a transformation of the
neighborhood of the point $(1,0, \ldots,0)$ to the one
of $(1, \ldots, 1)$.
For every element $u=(u_1, \ldots, u_g) \in (\xz / 2 \xz)^g$ we define a
multiplicative character
\[
l_u:(\xz / 2 \xz)^g \rightarrow \{ 1,-1 \}, \quad (v_1,\ldots,v_g)
\mapsto \prod_{i=1, \ldots, g} (-1)^{u_iv_i}.
\]
We note that $l_u(v)=l_v(u)$ for all $u,v \in (\xz / 2 \xz)^g$.
Now consider the linear transformation
which is given by the matrix
$M=(m_{u,v})_{u,v \in (\xz / 2 \xz)^g} \in \mathrm{Mat}(2^g,
\xz)$ with entries $m_{u,v}=l_u(v)$.
Obviously, the matrix $M$
maps the point $(1,0, \ldots , 0)$ to the point $(1, \ldots ,1)$.
\newline\indent
In the following we assume that we have normalized the sequence
of theta null points
$\big( a_u^{(n)} \big)$ such that equation (\ref{ctrelation}) holds with
$\lambda_n=1$ for all $n \geq 2$.
The following lemma was communicated to the author by
Lubicz. By lack of a suitable reference we give a proof in here.
\begin{lemma}
\label{trafo}
If we define $\big( x_u^{(n)} \big)_{u \in (\xz / 2 \xz)^g}$ and $\big(
x_u^{(n+1)} \big)_{u \in (\xz / 2 \xz)^g}$
coordinatewise by the formulas
\[
x_u^{(n)} = \left( \sum_{v \in (\xz / 2 \xz)^g}
m_{u,v} \cdot a_v^{(n)} \right)^2
\quad \mbox{and} \quad x_u^{(n+1)} = \left(
\sum_{v \in (\xz / 2 \xz)^g}
m_{u,v} \cdot a_v^{(n+1)}  \right)^2
\]
then $x_u^{(n)}$ and $x_u^{(n+1)}$ satisfy the
relation {\rm (\ref{gagm})} for all $u \in (\xz / 2 \xz)^g$.
\end{lemma}
\begin{proof}
We set
\[
y_u^{(n)} = \sum_{v \in (\xz / 2 \xz)^g}
m_{u,v} \cdot a_v^{(n)}
\quad \mbox{and} \quad y_u^{(n+1)} =
\sum_{v \in (\xz / 2 \xz)^g}
m_{u,v} \cdot a_v^{(n+1)}
\]
One computes
\begin{eqnarray}
\label{eeq1}
\lefteqn{\frac{1}{2^g} \cdot \sum_{v \in (\xz /2 \xz)^g} y_{u+v}^{(n)}
\cdot y_{v}^{(n)}} \\
\nonumber && = \frac{1}{2^g} \cdot \sum_{v \in (\xz /2 \xz)^g} \left(
\sum_{s \in (\xz /2 \xz)^g} l_{u+v}(s) \cdot a_s^{(n)} \right)
\cdot \left(
\sum_{t \in (\xz /2 \xz)^g} l_v(t) \cdot a_t^{(n)} \right)
\end{eqnarray}
\begin{eqnarray*}
\nonumber && = \frac{1}{2^g} \cdot \sum_{v \in (\xz /2 \xz)^g} \sum_{s,t \in (\xz /2 \xz)^g}
l_{u+v}(s) \cdot l_v(t) \cdot a_s^{(n)} \cdot a_t^{(n)} \\
\nonumber && = \frac{1}{2^g} \cdot \sum_{s,t \in (\xz /2 \xz)^g}
l_{u}(s) \cdot a_s^{(n)} \cdot a_t^{(n)}
\cdot \sum_{v \in (\xz /2 \xz)^g} l_v(s+t) \\
\nonumber && = \sum_{s \in (\xz /2 \xz)^g} l_{u}(s) \cdot \big( a_s^{(n)} \big)^2
\end{eqnarray*}
One the other hand we have
\begin{eqnarray}
\label{eeq2}
 \left( y_u^{(n+1)} \right)^2
=\left( \sum_{t \in (\xz /2 \xz)^g} l_u(t) \cdot a_t^{(n+1)}
\right)^2 = \sum_{s,t \in (\xz /2 \xz)^g} l_u(s+t) \cdot a_s^{(n+1)}
\cdot a_t^{(n+1)}
\end{eqnarray}
Let $w \in (\xz /2 \xz)^g$.
Applying equation (\ref{eeq1}) one gets
\begin{eqnarray}
\label{eeq3}
\lefteqn{\sum_{u \in (\xz /2 \xz)^g} l_w(u) \cdot
\left( \frac{1}{2^g} \cdot \sum_{v \in (\xz /2 \xz)^g} y_{u+v}^{(n)}
\cdot y_{v}^{(n)} \right)} \\
\nonumber && =
\sum_{u \in (\xz /2 \xz)^g} l_w(u) \cdot
\left( \sum_{s \in (\xz /2 \xz)^g} l_{u}(s) \cdot \big(
a_s^{(n)} \big)^2 \right)
=\sum_{s \in (\xz /2 \xz)^g} \big(
a_s^{(n)} \big)^2 \cdot \sum_{u \in (\xz /2 \xz)^g} l_w(u) \cdot l_u(s) \\
\nonumber &&
=\sum_{s \in (\xz /2 \xz)^g} \big(
a_s^{(n)} \big)^2 \cdot \sum_{u \in (\xz /2 \xz)^g} l_u(w+s)
= 2^g \cdot \big(
a_w^{(n)} \big)^2
\end{eqnarray}
Using equation (\ref{eeq2}) we compute
\begin{eqnarray}
\label{eeq4}
\lefteqn{\sum_{u \in (\xz /2 \xz)^g} l_w(u) \cdot
\left( y_u^{(n+1)} \right)^2} \\
\nonumber && = \sum_{u \in (\xz /2 \xz)^g} l_w(u) \cdot
\left( \sum_{s,t \in (\xz /2 \xz)^g} l_u(s+t) \cdot a_s^{(n+1)}
\cdot a_t^{(n+1)} \right) \\
\nonumber && =
\sum_{s,t \in (\xz /2 \xz)^g} a_s^{(n+1)}
\cdot a_t^{(n+1)} 
\sum_{u \in (\xz /2 \xz)^g} l_w(u) \cdot  l_u(s+t) \\
\nonumber && =
\sum_{s,t \in (\xz /2 \xz)^g} a_s^{(n+1)}
\cdot a_t^{(n+1)} 
\sum_{u \in (\xz /2 \xz)^g} l_u(w+s+t) \\
\nonumber && =
2^g \cdot \sum_{s \in (\xz /2 \xz)^g} a_s^{(n+1)}
\cdot a_{w+s}^{(n+1)}
\end{eqnarray}
By the equations (\ref{eeq3}) and (\ref{eeq4})
the equality (\ref{ctrelation}) holds for $\lambda_n=1$ if and only if
the relation (\ref{gagm}) is satisfied with the choice
of $x_u^{(n)}$ and $x_u^{(n+1)}$ as in the lemma.
\end{proof}
\noindent
The Lemma \ref{trafo} implies that, using the linear transformation $M$,
one can compute the sequence
$\big( a_u^{(n)} \big)_{n \geq 2}$ in terms of Mestre's
generalized AGM sequence.
For $n \geq 2$
let $x_u^{(n)}$ and $x_u^{(n+1)}$ be defined as in Lemma \ref{trafo}.
Then Theorem \ref{breuk} implies that the sequence
(\ref{fractionalagm}) converges
coordinatewise for $n \rightarrow \infty$.

\subsection{Perspectives}

We use the notation of the preceding section.
We note that the linear transformation which is induced by the matrix $M$
of Section \ref{agmscheme}
is not invertible over the ring $\xz_q$ because
of the following fact.
\begin{lemma}
We have $\mathrm{det}(M)=\pm 2^{(2^{g-1}g)}$.
\end{lemma}
\begin{proof}
We claim that $M^2=2^g \cdot I_{2^g}$,
where $I_{2^g}$ denotes the unit matrix of dimension $2^g$.
Let $b_{u,v}$ denote the entries of the Matrix $M^2$,
where $u,v \in (\xz / 2 \xz)^g$.
Then we have
\begin{eqnarray*}
\lefteqn{
b_{u,v}=\sum_{k=0}^g m_{u,k} \cdot m_{k,v} = \sum_{k=0}^g l_u(k) \cdot
l_k(v)} \\
&& = \sum_{k=0}^g l_{u+v}(k) = \left\{
\begin{array}{c@{, \quad }c}
2^g & u=v \\
0 & \mathrm{else}
\end{array}
\right.
\end{eqnarray*}
This proves the claim.
The lemma now follows from the the multiplicativity of the
determinant.
\end{proof}
\noindent
We conclude that the two theories, on one hand the
theory of canonical theta null points and on the other hand
Mestre's theory of the generalized AGM sequence, are not
equivalent.
The difference lies in the special fiber at the prime $2$.
\newline\indent
Now let $A$ be an ordinary abelian variety over a perfect field $k$
of characteristic $2$. Let $g=\mathrm{dim}(A)$.
Assume that we are given an ample symmetric line bundle $\pol$
of degree $1$ on $A$ and a symmetric theta structure $\Theta$
of type $(\xz / 2 \xz)^g$ for the line bundle $\pol^2$.
Let $(a_u)_{u \in (\xz / 2 \xz)^g}$ denote the theta null point
with respect to the triple $(A,\pol^2,\Theta)$.
\begin{lemma}
The theta null point $(a_u)$ equals the point $(1, 0, \ldots, 0)$.
\end{lemma}
\begin{proof}
We denote $Z_2=(\xz / 2 \xz)^g$.
Let $\hat{Z}_{2,k}$ denote the Cartier dual of the
finite constant group $Z_{2,k}$ associated to the
abstract group $Z_2$.
Recall that a theta structure is given by a
$\xg_{m,k}$-equivariant isomorphism
$\Theta:G ( Z_2) \stackrel{\sim}{\rightarrow} G(\pol^2)$,
where $G(Z_2)=\xg_{m,k} \times Z_{2,k} \times \hat{Z}_{2,k}$
is the standard Heisenberg group of type $Z_2$ and $G(\pol^2)$
the theta group of the line bundle $\pol^2$.
For the definition of the theta group and the standard Heisenberg
group see \cite[$\S$1]{mu66} and
Section \ref{thetalb}.
Let $G$ be the group of $k$-rational $\xg_{m,k}$-equivariant
automorphisms of the group $G(Z_n)$.
The group $G$ acts by composition on the set of
$k$-rational theta structures of type $Z_2$ for the line bundle
$\pol^2$.
Every element of the group $G$ is of the form
$(\tau,\delta)$, where $\tau:Z_{2,k} \times \hat{Z}_{2,k}
\rightarrow \xg_{m,k}$ is a character and
$\delta:Z_{2,k} \times \hat{Z}_{2,k}
\stackrel{\sim}{\rightarrow} Z_{2,k} \times \hat{Z}_{2,k}$
is an automorphism which is compatible with the
pairing
\[
e_2:(Z_{2,k} \times \hat{Z}_{2,k})^2 \rightarrow \xg_{m,k},
\quad \big( (x_1,l_1),(x_2,l_2) \big)
\mapsto \frac{l_2(x_1)}{l_1(x_2)}
\]
The given theta structure $\Theta$ induces a
trivialization
\begin{eqnarray}
\label{newtriv}
Z_{2,k} \stackrel{\sim}{\rightarrow} A[2]^{\mathrm{et}}
\end{eqnarray}
where $A[2]^{\mathrm{et}}$ denotes the maximal
\'{e}tale quotient of $A[2]$.
By \cite[Th.2.2]{ca07} there exists
a canonical theta structure $\Theta^{\mathrm{can}}$
of type $Z_2$ for the line bundle $\pol^2$ which depends
on the trivialization (\ref{newtriv}).
In the following we denote the canonical theta null
point with respect to the canonical theta
structure by $(a_u^{\mathrm{can}})_{u \in Z_2}$.
The action of $G$ on the set of $k$-rational
theta structures of type $Z_2$ for $\pol^2$
is faithful and transitive.
Hence there exists a unique $g=(\tau,\delta) \in G$
such that $g \star \Theta^{\mathrm{can}}$ equals the given
theta structure $\Theta$.
Since we are over a field of characteristic $2$, the
character $\tau$ has to be trivial.
We note that, because $k$ is assumed as perfect,
every finite flat group over $k$ splits
uniquely in a product of its connected component
with its maximal \'{e}tale quotient.
Hence the isomorphism $\delta$ acts
componentwise as a pair of automorphisms
$(\delta_1,\delta_2)$, where $\delta_2$ is the
inverse of the Cartier dual of $\delta_1$.
By the choice of the trivialization
(\ref{newtriv}) we conclude that
$\delta_1=\delta_2=\mathrm{id}$.
Thus we have $\Theta=\Theta^{\mathrm{can}}$.
Theorem \ref{2adictheta} implies that up to
normalization we have
\[
(a_u)=(a_u^{\mathrm{can}})=(1,0,\ldots,0).
\]
This completes the proof of the lemma.
\end{proof}
\noindent
This shows that the points of Mestre's AGM sequence
cannot be $2$-theta null points of abelian schemes since this is not
so modulo the prime $2$.
Recall that the points in Mestre's AGM sequence lie in
a $2$-adic neighborhood of the point $(1,\ldots,1)$.
We conjecture that Mestre's AGM sequence has
a moduli interpretation in terms of degenerate
abelian varieties, or, more strongly, in terms of semi-abelian schemes.
Our conjecture should follow
from a moduli interpretation of the `boundary' of Mumford's moduli space.
This `boundary' is canonically given by the points
in the projective
closure of the moduli space of abelian
varieties with theta structure under the canonical
projective embedding.

\subsection*{Leitfaden}

This article is structured as follows.
We recall some classically known facts about algebraic theta functions
in Section \ref{algtheta}.
In Section \ref{galoisprop} we prove
the Galois-theoretic property of the canonical theta
structure.
In Section \ref{main} we prove
theta relations which have as solutions the
canonical theta null points of canonical lifts
of ordinary abelian varieties over
a perfect field of characteristic $2$.
We give some examples of canonical theta null points in Section \ref{ccft}.

\subsection*{Notation}

Let $R$ be a ring, $X$ an $R$-scheme and $S$ an $R$-algebra.
By $X_{S}$ we denote the base extended scheme $X
\times_{\spec{R}} \spec{S}$.
Let $\mathcal{M}$ be a sheaf on $X$. Then we denote by $\mathcal{M}_S$
the sheaf 
that one gets by pulling back via the projection $X_S \rightarrow X$.
Let $\ab{I}{X}{Y}$ be a morphism of $R$-schemes. Then $I_S$ denotes the
morphism that is induced by $I$ via base extension with $S$.
\newline\indent
We use the same symbol for a scheme/module and
the associated fppf-sheaf.
By a \emph{group} we mean a group
object in the category of fppf-sheaves.
If a group $G$ is representable by a scheme and
the representing object has the property of
being  finite (flat, \'{e}tale, connected, etc.)
then we simply say that $G$ is a finite (flat, \'{e}tale,
connected, etc.) group.
Similarly we will say that
a morphism of groups is finite (faithfully flat, smooth, etc.)
if the groups are representable
and the induced morphism of schemes has the
corresponding property.
\newline\indent
A group (morphism of groups) is called finite locally free if
it is finite flat and of finite presentation.
The Cartier dual of a finite locally free commutative
group $G$ will be denoted by $G^D$.
The multiplication by an integer $n \in \xz$ on $G$ will be denoted by
$[n]$.
A finite locally free and surjective morphism between groups is called
an \emph{isogeny}.
By an \emph{elliptic curve} we mean an abelian scheme of relative
dimension $1$.
\newline\indent
If $G$ and $H$ are groups then we denote by $\underline{\mathrm{Hom}}(G,H)$
the fppf-sheaf which is defined by
$U \mapsto \mathrm{Hom} \big( G(U),H(U) \big)$, where the latter
denotes the group homomorphisms $G(U) \rightarrow H(U)$.

\section{Background on algebraic theta functions}
\label{algtheta}

In this section we recall some classically known
facts about algebraic theta functions.
We present Mumford's results on an algebraic theory of theta functions
in such a way that they apply to our situation.
The theory of algebraic theta functions was developed by
David Mumford in \cite{mu66}, \cite{mu67a} and \cite{mu67b}.

\subsection{Theta structures} 
\label{thetalb}

Let $A$ be an abelian scheme over a ring $R$ and $\pol$ a
line bundle on $A$.
Consider the morphism
\[
\varphi_{\pol}:A \rightarrow \mathrm{Pic}^0_{A/R},
\hsp x \mapsto \langle T_x^* \pol \otimes \pol^{-1} \rangle
\]
where $\langle \cdot \rangle$ denotes the class in $\mathrm{Pic}^0_{A/R}$.
We set $\check{A}=\mathrm{Pic}^0_{A/R}$. Note that $\check{A}$
is the dual of $A$ in the category of abelian schemes.
We denote the kernel of the morphism $\varphi_{\pol}$ by $A[\pol]$.
A line bundle $\pol$ on $A$ satisfies $A[\pol]=A$
if and only if its class is in $\mbox{Pic}^0_{A/R}(R)$.
Also it is well-known that if $\pol$ is relatively ample then
$\varphi_{\pol}$ is an isogeny. In the latter case we say that
$\pol$ has degree $d$ if $\varphi_{\pol}$ is fiberwise of
degree $d$.
Let $S$ be an $R$-algebra.
We define the theta group of $\pol$ as the functor
\[
G(\pol)(S)= \left\{ \hsp (x, \varphi) \hsp | \hsp x \in A[\pol](S),
\hsp \varphi:\pol_S \stackrel{\sim}{\rightarrow}
T^*_x \pol_S \hsp \right\}.
\]
The functor $G(\pol)$ has the structure of a group given by the group law
\[
\big( (y,\psi),(x,\varphi) \big) \mapsto (x+y, T_x^* \psi \circ \varphi).
\]
There are natural morphisms
\[
G(\pol) \rightarrow A[\pol], \hsp (x,\varphi) \mapsto x \quad \mbox{and}
\quad \xg_{m,R} \rightarrow G(\pol), \hsp \alpha \mapsto (0_A,\tau_{\alpha})
\]
where $0_A$ denotes the zero section of $A$ and $\tau_{\alpha}$ denotes
the automorphism of $\pol$ given by the multiplication with $\alpha$.
The induced sequence of groups
\begin{eqnarray}
\label{thetaseq}
0 \rightarrow \xg_{m,R} \rightarrow G(\pol) \stackrel{\pi}{\rightarrow} A[\pol]
\rightarrow 0
\end{eqnarray}
is central and exact.
Now let $\pol$ be relatively ample of degree $d$.
Then $A[\pol]$ is finite locally
free of order $d^2$.
\newline\indent
Let $K$ be a finite constant group.
The \emph{standard Heisenberg group} of type $K$ is given by the scheme
$\xg_{m,R} \times K_R \times K_R^D$ with group operation
\[
(\alpha_1,x_1,l_1) \star (\alpha_2,x_2,l_2) \stackrel{\mathrm{def}}{=}
(\alpha_1 \cdot \alpha_2 \cdot l_2(x_1), x_1+x_2,l_1 \cdot l_2).
\]
Here $K_R^D$ denotes the Cartier dual of $K_R$.
\begin{definition}
\label{deftheta}
A \emph{theta structure} of type $K$ for $\pol$ is a $\xg_{m,R}$-equivariant
isomorphism of groups
\[
G(K) \stackrel{\sim}{\rightarrow} G(\pol).
\]
\end{definition}
\noindent
For more details about theta structures
we refer to \cite[$\S$1]{mu66}.

\subsection{Mumford's Isogeny Theorem}
\label{isogtheo}

Let $R$ be a noetherian local ring.
Let $I:A \rightarrow A'$ be an isogeny of abelian schemes over $R$.
Assume that we are given ample line bundles $\pol$ and $\pol'$
on $A$ and $A'$, respectively.
Suppose we are given theta structures
\[
\Theta_A:G(K_{A}) \stackrel{\sim}{\rightarrow} G(\pol)
\quad \mbox{and} \quad
\Theta_{A'}:G(K_{A'}) \stackrel{\sim}{\rightarrow} G(\pol')
\]
where $K_A$ and $K_{A'}$ are finite constant groups.
Let $\alpha:I^* \pol' \stackrel{\sim}{\rightarrow} \pol$
be an isomorphism of $\mathcal{O}_A$-modules.
The existence of $\alpha$
implies that $\mathrm{Ker}(I)$ is contained in $A[\pol]$.
By \cite[$\S$1,Prop.1]{mu66} the morphism $\alpha$ induces
a section $\mathrm{Ker}(I) \rightarrow G(\pol)$
of the natural projection $G(\pol) \rightarrow A[\pol]$.
Let $\widetilde{K}$ denote the image of $\mathrm{Ker}(I)$ in $G(\pol)$.
Assume that
\begin{center}
\begin{tabular}{ll}
($\dagger$) & the image of $\widetilde{K}$ under $\Theta_A^{-1}$ is of the form $\{ 1 \} \times Z_1 \times Z_2$ \\
& with subgroups $Z_1 \leq K_A$ and $Z_2 \leq K_A^D$.
\end{tabular}
\end{center}
One defines
\[
Z_1^{\bot} = \{ \hsp x \in K_A \hsp | \hsp  (\forall l \in Z_2) \hsp l(x)=1 \hsp \}
\]
and
\[
Z_2^{\bot} = \{ \hsp l \in K_A^D \hsp | \hsp (\forall x \in Z_1) \hsp l(x)=1 \hsp \}.
\]
Note that $Z_1^{\bot} \times Z_2^{\bot}$ is the subgroup of points of
$K_A \times K_A^D$ being orthogonal to $Z_1 \times Z_2$ and that
$Z_1 \times Z_2 \subseteq Z_1^{\bot} \times Z_2^{\bot}$.
\begin{proposition}
\label{indtheta}
Let $\sigma:Z_1^{\bot} \rightarrow K_{A'}$ be
a surjective morphism of groups having kernel $Z_1$.
There exists a theta structure $\Theta_A(\sigma)$ of type $K_{A'}$ for the pair
$(A',\pol')$ depending on the theta structure $\Theta$ and the
morphism $\sigma$.
\end{proposition}
\begin{proof}
Let $G(\pol)^*$ denote the centralizer of $\widetilde{K}$ in $G(\pol)$.
There exists a natural isomorphism
of groups
\begin{eqnarray}
\label{mumiso}
G(\pol)^* / \widetilde{K} \stackrel{\sim}{\rightarrow} G(\pol')
\end{eqnarray}
(compare \cite[Prop.2]{mu66}).
The image of $G(\pol)^*$ under the morphism $\Theta_A^{-1}$ equals
$\xg_{m,R} \times Z_1^{\bot} \times Z_2^{\bot}$.
The isomorphism (\ref{mumiso}) composed with the induced isomorphism
\[
\xg_{m,R} \times Z_1^{\bot}/ Z_1 \times Z_2^{\bot} / Z_2
\stackrel{\sim}{\rightarrow} G(\pol)^*/ \tilde{K}
\]
gives an isomorphism
\begin{eqnarray}
\label{halftheta}
\xg_{m,R} \times Z_1^{\bot}/ Z_1
\times Z_2^{\bot} / Z_2\stackrel{\sim}{\rightarrow} G(\pol').
\end{eqnarray}
We claim that there exists a natural isomorphism
\begin{eqnarray}
\label{sigi}
Z_2^{\bot} / Z_2 \stackrel{\sim}{\rightarrow} (Z_1^{\bot}/ Z_1)^D.
\end{eqnarray}
This follows from the Snake Lemma applied to the following
commutative diagram of exact sequences
\[
\xymatrix{ 0 \ar@{->}[r] & Z_1 \ar@{->}[r] \ar@{->}[d]
& K_A \ar@{->}[d]^{\mathrm{id}} \ar@{->}[r] & (Z_2^{\bot})^D
\ar@{->}[d] \ar@{->}[r] & 0 \\
0 \ar@{->}[r] & Z_1^{\bot} \ar@{->}[r] & K_A \ar@{->}[r]^{\mathrm{e}} & Z_2^D
\ar@{->}[r] & 0. }
\]
Here the upper exact sequence is obtained by dualizing the exact
sequence
\[
0 \rightarrow Z_2^{\bot} \rightarrow K_A^D \stackrel{\mathrm{res}}{\longrightarrow}
Z_1^D \rightarrow 0
\]
and $\mathrm{e}$ is given by $x \mapsto \mathrm{e}
\big( (x,1),(0, \cdot ) \big)$ where $\mathrm{e}:(K_A \times K_A^D)^2
\rightarrow \xg_{m,R}$
denotes the commutator pairing.
The left hand vertical morphism is the natural inclusion.
The morphism $\sigma$ induces an isomorphism
$\sigma_1: K_{A'} \stackrel{\sim}{\rightarrow} Z_1^{\bot}/ Z_1$.
We denote by $\sigma_2: K_{A'}^D \stackrel{\sim}{\rightarrow} Z_2^{\bot}/ Z_2$
the inverse of the composition of
$\sigma_1^D$ with the isomorphism (\ref{sigi}).
Composing the isomorphism (\ref{halftheta}) with the isomorphism
$\mathrm{id} \times \sigma_1 \times \sigma_2$
we get a theta structure
\[
\Theta(\sigma):G(K_{A'}) \stackrel{\sim}{\rightarrow} G(\pol').
\]
This proves the proposition.
\end{proof}
\noindent
Note that $\Theta_A(\sigma)$ does not depend on the choice of $\alpha$.
\begin{definition}
\label{icompat}
We say that $\Theta_{A}$ and $\Theta_{A'}$ are $I$-compatible
if there exists $\alpha$ as above,
assumption {\rm ($\dagger$)} holds and
there exists a morphism $\sigma$ as above such that
$\Theta_{A'}=\Theta_A(\sigma)$.
\end{definition}
Let $\pi_A$ and $\pi_{A'}$ denote the structure maps of $A$ and $A'$,
respectively. 
Since $I$ is faithfully flat the natural morphism
$\pol' \stackrel{\tau}{\rightarrow} I_* I^* \pol'$ is injective.
As a consequence there exists an injective morphism $\iota:\pi_*' \pol'
\rightarrow \pi_* \pol$ of $\mathcal{O}_R$-modules
given by the composition
\[
\pi_*' \pol' \stackrel{\pi'_*\tau}{\longrightarrow}
\pi'_* I_* I^* \pol' = \pi_* I^* \pol' \stackrel{\pi_* \alpha}{\longrightarrow}
\pi_* \pol.
\]
The morphism $\iota$
identifies the sections of $\pi'_*\pol'$ with those sections of $\pi_* \pol$
which are invariant under the translations with points in the kernel of the
isogeny $I$.
Assume that we have chosen theta group invariant isomorphisms
\[
\beta_A:\pi_*\pol \stackrel{\sim}{\rightarrow} V(K_A)
\quad \mbox{and} \quad
\beta_{A'}:\pi_*'\pol'
\stackrel{\sim}{\rightarrow} V(K_{A'}),
\]
where we set $V(K)=\underline{\mathrm{Hom}}(K, \mathcal{O}_R)$ for 
a finite group $K$.
We define $V_I:V(K_{A'}) \rightarrow V(K_A)$
by setting
$V_I= \beta_A \circ \iota \circ \beta_{A'}^{-1}$.
\begin{theorem}[Isogeny Theorem]
\label{it}
Suppose $\Theta_A$ and $\Theta_{A'}$ are $I$-compatible.
In particular we are given a morphism $\sigma$
as above such that $\Theta_{A'}=\Theta_A(\sigma)$.
There exists a $\lambda \in R^*$ such that for all $f \in V(K_{A'})$
we have
\[
V_I(f)(x)=\left\{
\begin{array}{l@{, \quad}l}
0 & x \not\in Z_1^{\bot} \\
\lambda  f \big(\sigma(x) \big) & x \in Z_1^{\bot} 
\end{array}
\right.
\]
where $x \in K_A$.
\end{theorem}
\noindent
For a proof of Theorem \ref{it} in the case where $R$ is a field
see \cite[$\S$1, Th.4]{mu66}. The latter proof can easily be adapted
to our situation. 

\subsection{Mumford's $2$-multiplication formula}
\label{multform}

Let $R$ be a noetherian local ring and let $A \stackrel{\pi}{\rightarrow} \spec{R}$ be an abelian scheme
and $\pol$ an ample symmetric line bundle on $A$ which is
generated by global sections.
Let $n \geq 2$ be an integer. We set
$\pol_n=\pol^{\otimes n}$.
Assume that we are given
theta structures
\[
\Theta:G(K) \stackrel{\sim}{\rightarrow} G(\pol)
\quad \mbox{and} \quad
\Theta_n:G(K_n) \stackrel{\sim}{\rightarrow} G(\pol_n).
\]
One defines a morphism of groups
$\epsilon_n:G(\pol) \rightarrow G(\pol_n)$ by setting
\[
(x,\psi) \mapsto (x, \psi^{\otimes n}).
\]
Note that there is a natural inclusion $A[\pol] \hookrightarrow A[\pol_n]$
and the multiplication-by-$n$ induces an epimorphism $A[\pol_n] \rightarrow A[\pol]$.
On the image of $\xg_{m,R}$ the morphism
$\epsilon_n$ equals the $n$-th powering morphism.
\newline\indent
Next we define a morphism of groups
$\eta_n:G(\pol_n) \rightarrow G(\pol)$ using the symmetry of $\pol$.
Assume we are given $(x,\psi) \in G(\pol_n)$.
Since $\pol$ is symmetric there exists an isomorphism
$\gamma: \pol^{\otimes n^2} \stackrel{\sim}{\rightarrow} [n]^* \pol$.
Consider the composed isomorphism
\[
[n]^* \pol \stackrel{\gamma^{-1}}{\rightarrow} \pol^{\otimes n^2}
= \pol_n^{\otimes n} \stackrel{\psi^{\otimes n}}{\longrightarrow}
T_x^* \pol_n^{\otimes n} = T_x^* \pol^{\otimes n^2}
\stackrel{T^*_x \gamma}{\longrightarrow} T_x^* [n]^* \pol
=[n]^* T_{nx}^* \pol.
\]
Since $nx$ is a point of $A[\pol]$ there exists an isomorphism
$\rho: \pol \stackrel{\sim}{\rightarrow} T_{nx}^* \pol$
inducing the above isomorphism.
Since $[n]$ is faithfully flat the morphism
$\rho$ is uniquely determined.
We set $\eta_n(x, \psi)=(nx, \rho)$.
One can check that this definition is independent of the choice of $\gamma$.
The map $\eta_n$ restricted to $\xg_{m,R}$ equals the $n$-th powering
morphism.
\newline\indent
We denote the Lagrangian structures corresponding to $\Theta$ and $\Theta_n$
by $\delta$ and $\delta_n$, respectively.
Suppose that $K \leq K_n$ and $K=\left\{ nx| x \in K_n \right\}$.
Also we assume that $\delta_n$ restricted to $K \times K^D$
equals $\delta$.
As a consequence the multiplication-by-$n$ morphism on $K_n \times K_n^D$
induces an epimorphism
\[
K_n \times K_n^D \rightarrow K \times K^D,
(x,l) \mapsto (nx,l^n).
\]
We define morphisms
\[
E_n:G(K) \rightarrow G(K_n)
\quad \mbox{and} \quad
H_n:G(K_n) \rightarrow G(K)
\]
by setting
\[
(\alpha,x,l) \mapsto (\alpha^n,x,l)
\quad \mbox{and} \quad
(\alpha,x,l) \mapsto (\alpha^n,nx,l^n),
\]
respectively.
Here we consider points of $K \times K^D$ as points of $K_n \times K_n^D$
via the natural inclusion.
\begin{definition}
\label{compatibletheta}
We say that the theta structures $\Theta_n$ and $\Theta$ are $n$-compatible
if
\begin{enumerate}
\item[{\rm (i)}]
$K \leq K_n$ and $K=\left\{ nx| x \in K_n \right\}$,
\item[{\rm (ii)}]
$\delta_n$ restricted to $K \times K^D$
equals $\delta$,
\item[{\rm (iii)}]
$\Theta_n \circ E_n = \epsilon_n \circ \Theta$ and
$\Theta \circ H_n = \eta_n \circ \Theta_n$.
\end{enumerate}
\end{definition}
\noindent
Now assume that $n=2$.
Choose theta group invariant isomorphisms
\[
\beta : \pi_* \pol \stackrel{\sim}{\rightarrow} V(K)
\quad \mbox{and} \quad
\beta_2 : \pi_* \pol_2 \stackrel{\sim}{\rightarrow} V(K_2),
\]
where $V(K)$ and $V(K_2)$ are defined as in the preceding section.
\begin{definition}
Let $s$ and $s'$ be sections of $\pi_* \pol$.
Set $f=\beta(s)$ and $f'=\beta(s')$.
We define
\[
f \star f' = \beta_2(s \otimes s').
\]
\end{definition}
\noindent
Assume that we have chosen a rigidification of $\pol$.
The latter choice implies the existence of
finite theta functions $q_{\pol} \in V(K)$ and
$q_{\pol_2} \in V(K_2)$ which interpolate the coordinates
of the corresponding theta null points (cf.
\cite[$\S$1]{mu66}). 
\begin{theorem}[$2$-Multiplication Formula]
\label{mult}
Suppose $\Theta$ and $\Theta_2$ are $2$-compatible theta structures.
There exists a $\lambda \in R^*$ such that for all $f,f' \in V(K)$ we have
\[
(f \star f')(x)= \lambda \sum_{y \in x + K} f(x+y)
 f'(x-y)  q_{\pol_2}(y)
\]
where $x \in K_2$.
\end{theorem}
\noindent
For a proof of Theorem \ref{mult} over a field see \cite[$\S$3]{mu66}.
Without major changes the proof given there applies to
our situation.
\newline\indent
Assume that we have chosen an embedding
$Z_2 \subseteq K_2$
where $Z_2=(\ZZ / 2 \ZZ)^g$ and $g=\mathrm{dim}_R(A)$.
Applying the above theorem to the Dirac-basis of
$V(K)$ and evaluating at the zero section of the abelian scheme
$A$ gives the
following relation for theta null values.
\begin{corollary}
\label{qell2}
Suppose $\Theta$ and $\Theta_2$ are $2$-compatible theta structures.
Let $u,v \in K_2$ such that $u+v,u-v \in K$.
There exists a $\lambda \in R^*$ such that
\[
q_{\pol}(u+v) q_{\pol}(u-v) = \lambda  \sum_{z \in Z_2} q_{\pol_2}(u+z) q_{\pol_2}(v+z).
\]
\end{corollary}

\section{The Galois property of the canonical theta structure}
\label{galoisprop}

In this section we prove the characterizing pull back
property of the canonical theta structure (introduced in \cite{ca07}).
More precisely, we discuss the Galois theoretic properties of the
canonical theta structure in the
case where the theta structure is centered `at $p$', i.e. the
underlying level structure is non-\'{e}tale and the
corresponding polarization is inseparable.
For the \'{e}tale case see \cite[$\S$2.2]{ckl08}.

\subsection{Compatibility of the canonical theta structure}
\label{compatibility}

In the following we use the notation and the definitions of
Section \ref{algtheta}.
Let $R$ denote a complete noetherian local ring with perfect
residue field $k$ of characteristic $p>0$.
Assume that $R$ admits a lift of the $p$-th power Frobenius automorphism
of $k$.
Let $A \stackrel{\pi}{\rightarrow} \spec{R}$ be an abelian scheme
of relative dimension $g$
having ordinary reduction and
let $\pol$ be an ample symmetric line bundle of degree $1$ on $A$.
Let $F:A \rightarrow A^{(p)}$ denote the unique lift
of the relative $p$-Frobenius. We denote the structure map $A^{(p)}
\rightarrow \spec{R}$
by $\pi^{(p)}$.
By \cite[Th.5.1]{ca07} there exists an ample
symmetric line bundle $\polp$
of degree $1$ on $A^{(p)}$ and an isomorphism
$F^* \polp
\stackrel{\sim}{\rightarrow} \pol^{\otimes p}$
such that $(\polp)_k$ is the $p$-Frobenius twist of $\pol_k$.
For $i \geq 0$ we set
\[
\pol_i = \pol^{\otimes p^i}, \quad \emm_i= \big( \polp \big)^{ \otimes p^i}
\quad \mbox{and} \quad K_i=(\xz / p^i \xz)^g_R.
\]
Let $r \geq 1$ and assume that we are given an isomorphism
\begin{eqnarray}
\label{triv}
K_r \stackrel{\sim}{\rightarrow}
A[p^r]^{\mathrm{et}},
\end{eqnarray}
where $A[p^r]^{\mathrm{et}}$ denotes the maximal \'{e}tale quotient
of $A[p^r]$.
The lift of the relative $p$-Frobenius $F$ induces an
isomorphism $F[p^r]^{\mathrm{et}}:A[p^r]^{\mathrm{et}}
\stackrel{\sim}{\rightarrow} A^{(p)}[p^r]^{\mathrm{et}}$.
Composing $F[p^r]^{\mathrm{et}}$ with the isomorphism (\ref{triv})
we get an isomorphism
\begin{eqnarray}
\label{trivtwo}
K_r \stackrel{\sim}{\rightarrow} A^{(p)}[p^r]^{\mathrm{et}}.
\end{eqnarray}
Now assume that $A$ is the canonical lift of $A_k$.
As a consequence $A^{(p)}$ is also a canonical lift.
By \cite[Th.2.2]{ca07}
there exist for all $0 \leq j \leq r$
canonical theta structures $\Theta_j:G(K_j)
\stackrel{\sim}{\rightarrow} G(\pol_j)$ and
$\Sigma_j:G(K_j) \stackrel{\sim}{\rightarrow} G(\emm_j)$
depending on the
isomorphisms (\ref{triv}) and (\ref{trivtwo}), respectively.
For the definition of a theta structure see Definition \ref{deftheta}.
\newline\indent
Let us fix some notation that will be used in proofs of
the following lemmas.
By $V_i:A \rightarrow A_i$, where $i \geq 1$, we denote the $i$-fold application
of the lift of the $p$-Verschiebung.
For $1 \leq i \leq r$ let $\delta_i$ and $\kappa_i$ denote
the Lagrangian structures induced by the theta structures
$\Theta_i$ and $\Sigma_i$, respectively.
We note that by the definition of the canonical theta structure
we have $A_i=A/\delta_i(K_i)$ for $i \leq r$.
Let
\[
v_{j+1}:K_{j+1} \rightarrow G(\pol_{j+1})
\quad \mbox{and} \quad
v_j:K_j \rightarrow G(\pol_j)
\]
be the liftings induced by the canonical theta structures $\Theta_{j+1}$ and
$\Theta_j$, respectively.
By descent the sections $v_{j+1}$ and $v_j$
correspond to line bundles
$\pol^{(j+1)}$ and $\pol^{(j)}$ on $A_{j+1}$ and $A_j$,
respectively, equipped with isomorphisms
\begin{eqnarray*}
\beta_{j+1}:V_{j+1}^* \pol^{(j+1)} \stackrel{\sim}{\rightarrow}
\pol_{j+1}
\quad \mbox{and} \quad
\beta_j:V_j^* \pol^{(j)} \stackrel{\sim}{\rightarrow}
\pol_j.
\end{eqnarray*}
Also, we set
$A_0=A$ and $\pol^{(0)}= \pol$.
\begin{lemma}
\label{comp1}
For $0 \leq j < r$ the theta structures
$\Theta_{j+1}$ and $\Sigma_{j+1}$ are $p$-compatible
with $\Theta_j$ and $\Sigma_j$, respectively.
\end{lemma}
\begin{prf}{Proof of Lemma \ref{comp1}}
It suffices to prove the lemma for the theta structures
$\Theta_j$ where $0 \leq j \leq r$.
For trivial reasons the theta structure $\Theta_1$ is compatible with the
theta structure $\Theta_0$.
Assume that $j \geq 1$.
Obviously, the theta
structures $\Theta_{j+1}$ and $\Theta_j$
satisfy the conditions
(i) and (ii) of Definition \ref{compatibletheta}.
It remains to show that condition (iii) of the latter
definition is satisfied. We will proceed in two steps.
\begin{claim}
\label{claima}
One has
\begin{eqnarray}
\label{compeins}
\Theta_{j+1} \circ E_p = \epsilon_p \circ \Theta_j,
\end{eqnarray}
where $E_p$ and $\epsilon_p$ are defined as in
Section \ref{multform}.
\end{claim}
\begin{prf}{Proof of Claim \ref{claima}}
We verify equation (\ref{compeins}) for points which lie above $K_j$.
As the proof for points lying over $K_j^D$ is analogous we do not
present it here.
Let $x$ be a point of $K_j$.
One has
\[
v_{j+1}(x)= (\delta_{j+1}(x),T_{\delta_{j+1}(x)}^*\beta_{j+1} \circ \beta_{j+1}^{-1})
\quad \mbox{and} \quad
v_j(x)= (\delta_j(x),T_{\delta_j(x)}^*\beta_j \circ \beta_j^{-1}).
\]
Here we consider $x$ as a point of $K_{j+1}$ by the given inclusion
$K_j \subseteq K_{j+1}$.
By the definition of the canonical theta structure
there exists an isomorphism
\[
\beta:V^* \pol^{(j+1)}
\stackrel{\sim}{\rightarrow} (\pol^{(j)})^{\otimes p}.
\]
where $V:A_j \rightarrow A_{j+1}$ denotes the $p$-Verschiebung.
Consider the isomorphism $\psi$ given by the composition
\[
V_{j+1}^* \pol^{(j+1)} = V_j^*V^* \pol^{(j+1)}
\stackrel{V_j^* \beta}{\longrightarrow} V_j^* (\pol^{(j)})^{\otimes p}
\stackrel{\beta_j^{\otimes p}}
{\longrightarrow} \pol_j^{\otimes p} = \pol_{j+1}.
\]
Since $\delta_j(x)$ lies in the kernel of $V_j$
we have $V_j \circ T_{\delta_j(x)} = V_j$ and hence
$T_{\delta_j(x)}^*V_j^* \beta = V_j^* \beta$.
It follows that
\begin{eqnarray*}
& \epsilon_p ( v_j(x)) & = \epsilon_p \big( \delta_j(x),T_{\delta_j(x)}^*\beta_j \circ \beta_j^{-1}) \big) = \big( \delta_{j+1}(x),T_{\delta_{j+1}(x)}^*\beta_j^{\otimes p} \circ \beta_j^{ \otimes- p}
\big) \\
& & = \big( \delta_{j+1}(x),T_{\delta_{j+1}(x)}^*\psi \circ \psi^{-1} \big)
= \big(\delta_{j+1}(x),T_{\delta_{j+1}(x)}^*\beta_{j+1} \circ \beta_{j+1}^{-1} \big) = v_{j+1}(x).
\end{eqnarray*}
The latter equality sign follows from the fact that
$\beta_{j+1}$ and $\psi$ differ by a unit.
This proves our claim.
\end{prf}
\begin{claim}
\label{claimb}
One has
\begin{eqnarray}
\label{total}
\Theta_j \circ H_p = \eta_p \circ \Theta_{j+1}
\end{eqnarray}
where $H_p$ and $\eta_p$ are defined as in
Section \ref{multform}.
\end{claim}
\begin{prf}{Proof of Claim \ref{claimb}}
We verify the equality (\ref{total})
for points of $G(K_{j+1})$ which lie above $K_{j+1}$.
The proof for points lying over $K_{j+1}^D$ is analogous. 
Consider a point $(1,x,1)$ in $G(K_{j+1})$.
We have
\[
\Theta_j(H_p(1,x,1))=v_j(px)=(\delta_j(px),\tau_j)
\]
and
\[
\Theta_{j+1}(1,x,1)=v_{j+1}(x)=
(\delta_{j+1}(x),\tau_{j+1})
\]
where
\[
\tau_{j+1}=T_{\delta_{j+1}(x)}^*\beta_{j+1} \circ \beta_{j+1}^{-1}
\quad \mbox{and} \quad
\tau_j=T^*_{\delta_j(px)}\beta_j \circ \beta_j^{-1}.
\]
Choose an isomorphism $\gamma:\pol_j^{\otimes p^2} \stackrel{\sim}{\rightarrow}
[p]^* \pol_j$.
Consider the composed isomorphism
\begin{eqnarray*}
\lefteqn{[p]^* \pol_j \stackrel{\gamma^{-1}}{\longrightarrow}
\pol_j^{\otimes p^2}=\pol_{j+1}^{\otimes p} \stackrel{\tau_{j+1}^{\otimes p}}
{\longrightarrow}
(T_{\delta_{j+1}(x)}^*\pol_{j+1})^{\otimes p}} \\
& & =T_{\delta_{j+1}(x)}^* \pol_j^{\otimes p^2}
\stackrel{T_{\delta_{j+1}(x)}^* \gamma}{\longrightarrow} T_{\delta_{j+1}(x)}^*[p]^* \pol_j = [p]^* T_{\delta_j(px)}^* \pol_j.
\end{eqnarray*}
We claim that the latter isomorphism is induced by $\tau_j$.
By the definition of the canonical theta structure
there exists an isomorphism
\[
\xi:F^* \pol^{(j)} \stackrel{\sim}{\rightarrow}
(\pol^{(j+1)})^{\otimes p}
\]
where $F$ denotes the lift of the
relative $p$-Frobenius.
The composed isomorphism
\[
V_{j+1}^* (\pol^{(j+1)})^{\otimes p} \stackrel{V_{j+1}^* \xi^{-1}}{\longrightarrow}
V_{j+1}^* F^* \pol^{(j)}=V_j^* [p]^* \pol^{(j)}=
[p]^* V_j^* \pol^{(j)} \stackrel{[p]^* \beta_j}{\longrightarrow} [p]^* \pol_j.
\]
differs from $\gamma \circ \beta_{j+1}^{\otimes p}$ by a unit.
We conclude that
\begin{eqnarray*}
\lefteqn{T_{\delta_{j+1}(x)}^* \gamma \circ \tau_{j+1}^{\otimes p} \circ \gamma^{-1}
 = T_{\delta_{j+1}(x)}^*(\gamma \circ \beta_{j+1}^{\otimes p}) \circ
(\gamma \circ \beta_{j+1}^{\otimes p})^{-1}} \\
& & = T_{\delta_{j+1}(x)}^* \big( [p]^* \beta_j \circ V_{j+1}^* \xi^{-1} \big) \circ
\big( [p]^* \beta_j \circ V_{j+1}^* \xi^{-1} \big)^{-1}
= T_{\delta_{j+1}(x)}^*[p]^* \beta_j \circ [p]^* \beta_j^{-1} \\
& & = [p]^*T_{\delta_j(px)}^* \beta_j \circ [p]^* \beta_j^{-1}
= [p]^*(T_{\delta_j(px)}^*\beta_j \circ \beta_j^{-1}) = [p]^*\tau_j.
\end{eqnarray*}
Note that $T_{\delta_{j+1}(x)}^*V_{j+1}^* \xi^{-1} = V_{j+1}^* \xi^{-1}$ since $x$
is in the kernel of $V_{j+1}$.
\end{prf}
\noindent
This completes the proof of the lemma.
\end{prf}
\begin{lemma}
\label{comp2}
For all $0 \leq j < r$ the theta structures $\Theta_{j+1}$
and $\Sigma_j$ are $F$-compatible.
\end{lemma}
\begin{prf}{Proof of Lemma \ref{comp2}}
As the lemma is trivial for $j=0$ we may assume that $j>1$.
By assumption there exists an isomorphism
\[
\gamma_{j+1}:F^* \emm_j \stackrel{\sim}{\rightarrow} \pol_{j+1}.
\]
Obviously assumption ($\dagger$) of Section \ref{isogtheo}
holds with $Z_1=0$
and $Z_2=K_1^D$.
It follows that $Z_2^{\bot}=K_{j+1}^D$.
By duality we conclude that $Z_1^{\bot}$ coincides with the
image of $K_j$ in $K_{j+1}$ under the given inclusion.
\begin{claim}
\label{claimc}
One has
\[
\Sigma_j=\Theta_{j+1}(\mathrm{id})
\]
where the notation is as in
Proposition \ref{indtheta}.
\end{claim}
\begin{prf}{Proof of Claim \ref{claimc}}
Checking the claim amounts to prove the commutativity of the diagram
\begin{eqnarray*}
\xymatrix{
G( \pol_{j+1} ) / \tilde{K} \ar@{<-}[rr]
\ar@{->}[d]^{\mathrm{can}} & & \xg_{m,R} \times Z_1^{\bot} / Z_1
\times Z_2^{\bot} / Z_2 \ar@{<-}[d]^{\mathrm{id} \times \sigma_1
\times \sigma_2} \\
G(\emm_j) \ar@{<-}[rr]^{\Sigma_j} & & \xg_{m,R} \times K_j \times K_j^D }
\end{eqnarray*}
where the upper horizontal arrow is induced by $\Theta_{j+1}$,
the morphism $\sigma_1$ equals the `identity' and
the morphism $\sigma_2$ equals its dual (compare proof of Proposition
\ref{indtheta}).
In the following we verify the commutativity of the above diagram
for points of the form $(1,x,1)\nolinebreak\in\nolinebreak\xg_{m,R} \times K_j \times K_j^D$.
An analogous proof exists for points of the form $(1,0,l)$.
Via the morphisms $\sigma_1$ and $\sigma_2$
we can consider $(1,x,1)$ as a point of
$G(K_{j+1})$.
Here we consider $x$ as an element of $K_{j+1}$ via the given
inclusion $K_j \subseteq K_{j+1}$.
By definition its image under $\Theta_{j+1}$ is given by
$v(x)$ where
\[
v:K_{j+1} \rightarrow G(\pol_{j+1})
\]
is the lifting which is induced by the existence of an isomorphism
\[
\alpha_{j+1}:V_{j+1}^* \pol^{(j+1)}
\stackrel{\sim}{\rightarrow} \pol_{j+1}.
\]
The class of $v(x)$ in $G( \pol_{j+1} ) / \tilde{K}$ can be
represented by an element of the form
\[
\big( \delta_{j+1}(y),T^*_{\delta_{j+1}(y)} \alpha_{j+1} \circ \alpha_{j+1}^{-1})
\]
where $y \in K_{j+1}$ is chosen such that $F \big( \delta_{j+1}(y) \big)=
\kappa_{j+1}(x)$.
On the other hand we have
$\Sigma_j(1,x,1)=w(x)$ where
$w:K_j \rightarrow G(\emm_j)$
is the lifting that is induced by the canonical theta structure.
The lifting $w$ corresponds to the line bundle $\pol^{(j-1)}$ on $A_{j-1}$
with isomorphism
\[
\beta_j:V_j^*\pol^{(j-1)} \stackrel{\sim}{\rightarrow} \emm_j.
\]
We have
\[
w(x)=(\kappa_{j}(x),T^*_{\kappa_{j}(x)}\beta_j \circ \beta_j^{-1}).
\]
By the definition of the canonical theta structure there exist isomorphisms
\[
\xi_1:F^* \pol^{(j-1)} \stackrel{\sim}{\rightarrow} (\pol^{(j)})^{\otimes p}
\quad \mbox{and} \quad \xi_2:V^* \pol^{(j+1)} \stackrel{\sim}{\rightarrow}
(\pol^{(j)})^{\otimes p}.
\]
Let $\xi=\xi_1^{-1} \circ \xi_2$. The isomorphism $V_j^* \xi$
induces an isomorphism
\[
V_{j+1} \pol^{(j+1)} = V_j^* V^* \pol^{(j+1)} \stackrel{V_j^* \xi}
{\longrightarrow} V_j^* F^* \pol^{(j-1)} = F^* V_j^* \pol^{(j-1)}.
\]
The composed isomorphism $\gamma_{j+1} \circ F^* \beta_j \circ V^*_j \xi$
differs from $\alpha_{j+1}$ by a unit.
By the definition of the canonical isomorphism
\[
G(\pol_{j+1}) / \tilde{K}
\stackrel{\sim}{\rightarrow} G(\emm_j)
\]
(see proof of \cite[$\S$1,Prop.2]{mu66})
the element in $G(\pol_{j+1})/\tilde{K}$ that corresponds to $w(x)$
is given by
\begin{eqnarray*}
\lefteqn{ \big( \delta_{j+1}(y),T^*_{\delta_{j+1}(y)} \gamma_{j+1} \circ F^* (T_{\kappa_{j}(x)}^*\beta_j
\circ \beta_j^{-1}) \circ \gamma_{j+1}^{-1} \big)} \\
& & = \big( \delta_{j+1}(y),T^*_{\delta_{j+1}(y)}\gamma_{j+1} \circ F^*T_{\kappa_j(x)}^*\beta_j
\circ F^*\beta_j^{-1} \circ \gamma_{j+1}^{-1} \big) \\
& & = \big( \delta_{j+1}(y),T^*_{\delta_{j+1}(y)}\gamma_{j+1} \circ T_{\delta_{j+1}(y)}^*F^*\beta_j
\circ F^*\beta_j^{-1} \circ \gamma_{j+1}^{-1} \big) \\
& & = \big(\delta_{j+1}(y),T^*_{\delta_{j+1}(y)} (\gamma_{j+1} \circ F^*\beta_j
\circ V_j^* \xi)
\circ (\gamma_{j+1} \circ F^*\beta_j \circ V_j^* \xi)^{-1} \big) \\
& & = (\delta_{j+1}(y),T^*_{\delta_{j+1}(y)}\alpha_{j+1} \circ \alpha_{j+1}^{-1}) = v(x).
\end{eqnarray*}
Hence our claim is proved.
\end{prf}
This completes the proof of the lemma.
\end{prf}

\subsection{The pull back of the canonical theta null point}
\label{pback}

Let $R$ be a complete noetherian local ring with
perfect residue field $k$ of characteristic $p>0$.
The following lemma forms the key ingredient of the
theorem that is proved below.
\begin{lemma}
\label{homnull}
One has
\[
\underline{\mathrm{Hom}} \big( (\QQ_p / \ZZ_p)_R^g,
\mu^g_{p^{\infty},R} \big)(R)=0.
\]
\end{lemma}
\begin{proof}
It suffices to prove that
\begin{eqnarray}
\label{nohom}
\underline{\mathrm{Hom}} \big( (\QQ_p / \ZZ_p)_R,
\mu_{p^{\infty},R} \big)(R)=0.
\end{eqnarray}
The elements of the group $\underline{\mathrm{Hom}} \big( (\QQ_p / \ZZ_p)_R,
\mu_{p^{\infty},R} \big)(R)$ correspond to the compatible systems $(a_1,a_2,a_3,\ldots)$ of $p$-power roots of unity in $R$,
i.e. $a_{i+1}^p=a_i$ for $i \geq 1$.
Such a compatible system necessarily has to be trivial. The latter is
obvious on the level Artin local quotients of $R$ by powers of the maximal
ideal.
This completes the proof of the remark.
\end{proof}
\noindent
Suppose that we are given an abelian scheme $A$ over $R$
which has ordinary reduction.
Let $\pol$ be an ample symmetric line bundle of degree $1$ on $A$.
We set $q=p^d$ where $d \geq 1$ is an integer.
Assume that there exists a $\sigma \in \mathrm{Aut}(R)$ lifting the $p$-th power
Frobenius automorphism of $k$.
Recall that there exists a canonical lift $F:A \rightarrow A^{(p)}$ of
the relative $p$-Frobenius morphism and a canonical ample symmetric line bundle
$\pol^{(p)}$ of degree $1$ on $A^{(p)}$ such that $F^* \pol^{(p)} \cong \pol^{p}$
(see \cite[$\S$5]{ca07}).
We set $Z_n = (\xz / n \xz)^g$ for $n \geq 1$ where $g=\mathrm{dim}_R(A)$.
Assume that we are given an isomorphism
\begin{eqnarray}
\label{trivpq}
Z_{q,R} \stackrel{\sim}{\rightarrow} A[q]^{\mathrm{et}},
\end{eqnarray}
where $A[q]^{\mathrm{et}}$ denotes the maximal \'{e}tale quotient of $A[q]$.
The Frobenius lift $F$ induces an isomorphism
\begin{eqnarray*}
F[q]^{\mathrm{et}}:
A[q]^{\mathrm{et}} \stackrel{\sim}{\rightarrow} A^{(p)}[q]^{\mathrm{et}}. 
\end{eqnarray*}
Composing $F[q]^{\mathrm{et}}$ with the isomorphism (\ref{trivpq})
gives an isomorphism
\begin{eqnarray}
\label{trivpq2}
Z_{q,R} \stackrel{\sim}{\rightarrow} A^{(p)}[q]^{\mathrm{et}}.
\end{eqnarray}
Assume that $A$ is the canonical lift of $A_k$.
Our assumption implies that the abelian scheme
$A^{(p)}$ is a canonical lift.
By \cite[Th.2.2]{ca07} there exist canonical theta structures
$\Theta_q$ and $\Theta^{(p)}_q$ of type $Z_q$
for the line bundles $\pol^q$ and $\big( \pol^{(p)} \big)^q$, which
depend on the isomorphisms (\ref{trivpq}) and (\ref{trivpq2}),
respectively.
\newline\indent
Let $A^{(\sigma)}$ be defined by the Cartesian diagram
\begin{eqnarray}
\label{defpullback}
\xymatrix{
A^{(\sigma)} \ar@{->}[r]^{\mathrm{pr}} \ar@{->}[d] & A \ar@{->}[d] \\
\spec{R} \ar@{->}[r]^{\spec{\sigma}} & \spec{R},}
\end{eqnarray}
where the right hand vertical arrow is the structure morphism.
Let $\pol^{(\sigma)}$ be the pull back of $\pol$ along the morphism
$\mathrm{pr}:A^{(\sigma)} \rightarrow A$ which is defined by
the diagram (\ref{defpullback}).
We obtain a theta structure $\Theta_q^{(\sigma)}$ for $\big( \pol^{(\sigma)} \big)^q$ by extension of scalars along $\spec{\sigma}$
applied to the theta structure $\Theta_q:G(Z_q) \stackrel{\sim}{\rightarrow}
G(\pol^q)$ and by chaining with the natural
isomorphism $Z_q \stackrel{\sim}{\rightarrow} Z_{q,\sigma}$ times
the inverse of its dual.
In the following we will assume
that the special fibers of $A^{(p)}$
and $A^{(\sigma)}$ are in fact equal.
By uniqueness of the canonical lift there exists a canonical isomorphism
$\tau: A^{(p)} \stackrel{\sim}{\rightarrow} A^{(\sigma)}$ over $R$
lifting the identity on special fibers.
We set $\pol_n= \pol^n$, $\pol^{(p)}_n = \big( \pol^{(p)} \big)^n$
and $\pol^{(\sigma)}_n = \big( \pol^{(\sigma)} \big)^n$ for $n \geq 1$. 
\begin{lemma}
We have $\tau^* \pol_q^{(\sigma)} \cong \pol_q^{(p)}$.
\end{lemma}
\begin{proof}
We set $\mathcal{M}= \big( \tau^* \pol^{(\sigma)} \big) \otimes \big(
\pol^{(p)} \big)^{-1}$.
It follows by the definition of $\pol^{(p)}$
that the class of $\mathcal{M}$ reduces
to the trivial class.
By assumption the line bundle $\tau^* \pol^{(\sigma)}$ is symmetric.
By \cite[Th.5.1]{ca07} also the line bundle $\pol^{(p)}$ is symmetric.
As a consequence, the line bundle $\mathcal{M}$ is symmetric
and gives an element of
$\mathrm{Pic}^0_{A^{(p)}/R}[2](R)$.
If $p=2$ then clearly $\mathcal{M}^q$ gives the trivial class.
Suppose that $p>2$.
We observe that by assumption the group $\mathrm{Pic}^0_{A^{(p)}/R}[2]$ is finite \'{e}tale.
We conclude by the connectedness of the ring $R$ and by the fact that
$\mathcal{M}$ gives the trivial class modulo $p$
that the class of $\mathcal{M}$ is trivial. This proves the lemma.
\end{proof}
\noindent
By the above discussion there exists an isomorphism
$\gamma: \tau^* \pol_q^{(\sigma)} \stackrel{\sim}{\rightarrow} \pol_q^{(p)}$.
We define a $\xg_{m,R}$-invariant morphism of theta groups $\tau^{\dagger}:G \big( \pol_q^{(\sigma)} \big) \rightarrow
G \big( \pol_q^{(p)} \big)$ by setting
\[
(x, \varphi) \mapsto \big( y, T_y^* \gamma \circ \tau^* \varphi \circ \gamma^{-1} \big)
\]
where $y = \tau^{-1}(x)$ and $\varphi: \pol_q^{(\sigma)}
\stackrel{\sim}{\rightarrow} T_x  \pol_q^{(\sigma)}$ is a given isomorphism.
Obviously, our definition is independent of the choice of
$\gamma$. For trivial reasons the morphism $\tau^{\dagger}$ gives an
isomorphism. By what has been said before, the isomorphism
$\tau^{\dagger}$ is canonical.
\begin{theorem}
\label{cantwisttheta}
One has
\begin{eqnarray}
\label{pbprop}
\tau^{\dagger} \circ \Theta_q^{(\sigma)} = \Theta_q^{(p)}.
\end{eqnarray}
\end{theorem}
\begin{prf}{Proof of Theorem \ref{cantwisttheta}}
In the following we fix the notation that we will use during the
course of the proof.
Let $\delta_q$, $\delta_q^{(p)}$ and $\delta_q^{(\sigma)}$ be the Lagrangian
structures which are induced by $\Theta_q$, $\Theta_q^{(p)}$ and
$\Theta_q^{(\sigma)}$, respectively.
The isomorphism $\tau$ induces isomorphisms
$\tau_q:A^{(p)}[q] \stackrel{\sim}{\rightarrow} A^{(\sigma)}[q]$
and $\tau_q^{\mathrm{et}}:A^{(p)}[q]^{\mathrm{et}} \stackrel{\sim}{\rightarrow}
A^{(\sigma)}[q]^{\mathrm{et}}$, where $A^{(p)}[q]^{\mathrm{et}}$
and $A^{(\sigma)}[q]^{\mathrm{et}}$ denote the maximal \'{e}tale
quotients of $A^{(p)}[q]$ and $A^{(\sigma)}[q]$, respectively.
Let
\[
\pi_q:A[q] \rightarrow A[q]^{\mathrm{et}}, \quad
\pi^{(\sigma)}_q: A^{(\sigma)}[q] \rightarrow A^{(\sigma)}[q]^{\mathrm{et}}
\quad \mbox{and} \quad
\pi^{(p)}_q: A^{(p)}[q] \rightarrow A^{(p)}[q]^{\mathrm{et}}
\]
denote the natural projections on maximal \'{e}tale quotients.
We denote the natural inclusion $Z_{q,R} \rightarrow Z_{q,R} \times Z_{q,R}^D$
by $i$.
\newline\indent
Furthermore, let $r_q^{(p)}:A^{(p)}[q]^{\mathrm{et}} \rightarrow A^{(p)}[q]$
and $r_q^{(\sigma)}:A^{(\sigma)}[q]^{\mathrm{et}} \rightarrow A^{(\sigma)}[q]$
be the sections of the natural projections $\pi_q^{(p)}$ and $\pi_q^{(\sigma)}$,
corresponding to the theta structures
$\Theta_q^{(p)}$ and $\Theta_q^{(\sigma)}$, respectively,
such that
\[
r_q^{(p)} \circ \pi^{(p)}_q \circ \delta_q^{(p)} \circ i
= \delta_q^{(p)} \circ i
\]
and
\[
r_q^{(\sigma)} \circ \pi^{(\sigma)}_q \circ \delta_q^{(\sigma)} \circ i
= \delta_q^{(\sigma)} \circ i.
\]
The proof of the theorem is divided into several lemmas.
We start the proof by checking that the equality (\ref{pbprop}) holds
for the induced Lagrangian structures.
\begin{lemma}
\label{lemmaa}
One has the equality
\begin{eqnarray}
\label{pullbacklag}
\delta_q^{(\sigma)} = \tau_q \circ \delta_q^{(p)}.
\end{eqnarray}
\end{lemma}
\begin{prf}{Proof of Lemma \ref{lemmaa}}
First we verify the equality (\ref{pullbacklag}) on points
which lie in the image of the morphism $i$.
Thus we have to show that
\begin{eqnarray}
\label{ieq}
\delta_q^{(\sigma)} \circ i= \tau_q \circ \delta_q^{(p)} \circ i.
\end{eqnarray}
\begin{claim}
\label{cll}
There is an equality
\begin{eqnarray}
\label{etpart}
\pi^{(\sigma)}_q \circ \delta_q^{(\sigma)} \circ i = \pi^{(\sigma)}_q \circ \tau_q \circ \delta_q^{(p)} \circ i
=  \tau_q^{\mathrm{et}} \circ \pi^{(p)}_q \circ \delta_q^{(p)} \circ i.
\end{eqnarray}
\end{claim}
\begin{prf}{Proof of Claim \ref{cll}}
The right hand equality of the equation (\ref{etpart}) is obvious.
In the following we will show the left hand equality of (\ref{etpart}).
Let $F_q^{\mathrm{et}}:A[q]^{\mathrm{et}} \rightarrow A^{(p)}[q]^{\mathrm{et}}$
denote the isomorphism which is induced by the Frobenius lift $F$.
By the definition of $\delta_q^{(p)}$ we have
\[
\pi^{(p)}_q \circ \delta_q^{(p)} \circ i =  F_q^{\mathrm{et}} \circ \pi_q  \circ \delta_q \circ i,
\]
and hence, by means of the right hand equality of (\ref{etpart}), the left hand equation of (\ref{etpart}) is equivalent to
\begin{eqnarray}
\label{tttheta}
\pi^{(\sigma)}_q \circ \delta_q^{(\sigma)} \circ i =  \tau_q^{\mathrm{et}} \circ
F_q^{\mathrm{et}} \circ \pi_q  \circ \delta_q \circ i.
\end{eqnarray}
By the theory of finite \'{e}tale
groups it suffices to check that the latter equality holds
on the special fiber.
We claim that this indeed is the case.
Our claim follows from the fact that by assumption
the morphism $\tau_q^{\mathrm{et}}$ equals
the identity on the special fiber.
\end{prf}
\begin{claim}
\label{clll}
One has
\begin{eqnarray}
\label{sectc}
\tau_q \circ r_q^{(p)} = r_q^{(\sigma)} \circ \tau_q^{\mathrm{et}}.
\end{eqnarray}
\end{claim}
\begin{prf}{Proof of Claim \ref{clll}}
The sections $r_q^{(\sigma)}$ and
$r_q^{(p)}$ are, by the definition of the canonical theta structure,
the truncations to $A^{(\sigma)}[q]^{\mathrm{et}}$
and $A^{(p)}[q]^{\mathrm{et}}$ of sections
$r_{p^{\infty}}^{(\sigma)}$ and $r_{p^{\infty}}^{(p)}$, respectively,
of the connected-\'{e}tale sequences of Barsotti-Tate groups
\begin{eqnarray}
\label{esone}
\xymatrix{
0 \ar@{->}[r] & A^{(\sigma)}[p^{\infty}]^{\mathrm{loc}}
\ar@{->}[r] & A^{(\sigma)}[p^{\infty}] \ar@{->}[r] &
A^{(\sigma)}[p^{\infty}]^{\mathrm{et}}
\ar@{->}[r] \ar@/_1pc/[l]_{r_{p^{\infty}}^{(\sigma)}} & 0}
\end{eqnarray}
and
\begin{eqnarray}
\label{estwo}
\xymatrix{
0 \ar@{->}[r] & A^{(p)}[p^{\infty}]^{\mathrm{loc}}
\ar@{->}[r] & A^{(p)}[p^{\infty}] \ar@{->}[r] &
A^{(p)}[p^{\infty}]^{\mathrm{et}}
\ar@{->}[r] \ar@/_1pc/[l]_{r_{p^{\infty}}^{(p)}} & 0.}
\end{eqnarray}
Proving the equality (\ref{sectc}) amounts to prove
the uniqueness
of the sections
$r_{p^{\infty}}^{(\sigma)}$ and $r_{p^{\infty}}^{(p)}$.
The uniqueness follows from the following observation.
Passing to the strict Henselization of $R$ we have
\[
A^{(\sigma)}[p^{\infty}]^{\mathrm{loc}} \cong
A^{(p)}[p^{\infty}]^{\mathrm{loc}} \cong \mu^g_{p^{\infty},R}
\quad \mbox{and} \quad
A^{(\sigma)}[p^{\infty}]^{\mathrm{et}} \cong
A^{(p)}[p^{\infty}]^{\mathrm{et}} \cong (\QQ_p / \ZZ_p)_R^g,
\]
where $g$ denotes the relative dimension of $A$ over $R$.
The sections of the exact sequence (\ref{esone}) and
(\ref{estwo}) form a
$\underline{\mathrm{Hom}} \big( (\QQ_p / \ZZ_p)_R^g,
\mu^g_{p^{\infty},R} \big)$-torsor.
Hence the uniqueness follows from Lemma \ref{homnull}.
\end{prf}
\noindent
Combining the Claims \ref{cll} and \ref{clll}
we conclude that the equality (\ref{ieq}) holds.
We omit a detailed proof of the equality (\ref{pullbacklag})
on the points of the connected component and refer to \cite[$\S$6.3]{ca07}
for details of the construction of the canonical level structure.
The proof essentially is a consequence of the fact that
the `connected factor' of the canonical level structure
equals (up to canonical isomorphism) the dual of the `\'{e}tale factor' of the
canonical level structure.
This completes the proof of the lemma.
\end{prf}
\noindent
It remains to show, on top of Lagrangian structures, that the equality of theta structures holds, as claimed in the theorem.
By \cite[Prop.4.2]{ca07} and \cite[Prop.4.5]{ca07} the verification of the equality
of the theta structures $\tau^{\dagger} \circ \Theta_q^{(\sigma)}$ and
$\Theta_q^{(p)}$ comes down to checking that certain
descent line bundles are isomorphic.
Let $F_q:A^{(p)}\rightarrow A^{(pq)}$ and $V_q:A^{(p)} \rightarrow
A^{(p/q)}$ be the $q$-fold application of the lifts
of the relative $p$-Frobenius and the Verschiebung, respectively.
By construction we have
\[
\mathrm{Ker}(F_q)=\delta^{(p)} \big( \{ 0 \} \times Z_{q,R}^D \big)
\quad \mbox{and} \quad 
\mathrm{Ker}(V_q)=\delta^{(p)} \big( Z_{q,R} \times \{ 0 \} \big).
\]
We first consider the descent of $\big( \pol^{(p)} \big)^q$
along $F_q$.
The theta structures $\tau^{\dagger} \circ \Theta_q^{(\sigma)}$
and $\Theta_q^{(p)}$ correspond to line bundles
$\mathcal{M}^{(\sigma)}$ and $\mathcal{M}^{(p)}$ on
$A^{(pq)}$ such that
\[
\big( \pol^{(p)} \big)^q \cong F_q^* \mathcal{M}^{(\sigma)}
\cong F_q^* \mathcal{M}^{(p)}.
\]
By definition of the canonical theta structure both of the line bundles
$\mathcal{M}^{(\sigma)}$ and $\mathcal{M}^{(p)}$
reduce to the $q$-Galois twist of $\pol_k^{(p)}$.
It follows by uniqueness \cite[Th.5.1]{ca07}
that $\mathcal{M}^{(\sigma)} \cong \mathcal{M}^{(p)}$.
\newline\indent
Now consider the case of descent along Verschiebung $V_q$.
The symmetric theta structures $\Theta_q^{(p)}$ and
$\tau^{\dagger} \circ \Theta_q^{(\sigma)}$
correspond to symmetric line bundles
$\mathcal{H}^{(p)}$ and $\mathcal{H}^{(\sigma)}$
of degree $1$, respectively, such that
\[
\big( \pol^{(p)} \big)^q \cong V_q^* \mathcal{H}^{(p)}
\cong V_q^* \mathcal{H}^{(\sigma)}.
\]
We claim that $\mathcal{H}^{(p)} \cong \mathcal{H}^{(\sigma)}$.
If $p>2$ then by symmetry
$\mathcal{H}^{(p)} \otimes \big( \mathcal{H}^{(\sigma)} \big)^{-1}$
gives an element of $\mathrm{Pic}^0_{A^{(q/p)}/R}[2](R)$ which lies in
the kernel of the dual $V_q^*$ of $q$-Verschiebung $V_q$. The kernel
of $V_q^*$ has
$p$-power order.
Because of the assumption $(p,2)=1$ the claim follows.
More delicate is the case where $p=2$.
Recall that the line bundles $\mathcal{H}^{(\sigma)}$ and $\mathcal{H}^{(p)}$
are uniquely determined by the choice of sections $s^{(\sigma)}$ and $s^{(p)}$
of the connected-\'{e}tale sequence
\begin{eqnarray*}
\xymatrix{
0 \ar@{->}[r] & A^{(2)}[2]^{\mathrm{loc}} \ar@{->}[r] &
A^{(2)}[2] \ar@{->}[r] & A^{(2)}[2]^{\mathrm{et}}
\ar@/_1pc/[l]_{s^{(\sigma)}} \ar@/^{1pc}/[l]^{s^{(p)}} \ar@{->}[r] & 0.
}
\end{eqnarray*}
(cf. \cite[Th.5.2]{ca07}).
The sections $s^{(p)}$ and $s^{(\sigma)}$,
which correspond to the theta
structures $\Theta_q^{(p)}$ and $\tau^{\dagger} \circ \Theta_q^{(\sigma)}$,
by construction both can be prolonged to sections of the
connected-\'{e}tale sequence
\begin{eqnarray*}
0 \rightarrow A^{(2)}[2^{\infty}]^{\mathrm{loc}}
\rightarrow A^{(2)}[2^{\infty}]
\rightarrow A^{(2)}[2^{\infty}]^{\mathrm{et}}
\rightarrow 0.
\end{eqnarray*}
But there exists a unique section of the latter exact sequence
of Barsotti-Tate groups.
By Lemma \ref{homnull} we conclude that
$s^{(p)}=s^{(\sigma)}$.
This completes the proof of the theorem.
\end{prf}

\section{Theta null points of $2$-adic canonical lifts}
\label{main}

Let $R$ be a complete noetherian local ring with perfect
residue field $k$ of characteristic $2$.
Assume that $R$ admits a lift $\sigma$
of the $2$nd power Frobenius automorphism
of $k$.
Let $A$ be an abelian scheme over $R$ of relative dimension $g$.
Let $\pol$ be an ample symmetric line bundle of degree $1$ on $A$.
We set $Z_n= (\ZZ / n \ZZ)^g$ for $n \geq 1$.
Assume that we are given $q=2^j$ where $j \geq 1$ and an isomorphism
\begin{eqnarray}
\label{trivia}
Z_{q,R} \stackrel{\sim}{\rightarrow}
A[q]^{\mathrm{et}}
\end{eqnarray}
where $A[q]^{\mathrm{et}}$ denotes the maximal \'{e}tale quotient
of $A[q]$.
Suppose that $A$ is the canonical lift of $A_k$.
By \cite[Th.2.2]{ca07} there exists a canonical theta
structure $\Theta_{q}$ of type $Z_{q}$ for
the pair $(A, \pol^{q})$ depending on the
isomorphism (\ref{trivia}).
Let $(x_u)_{u \in Z_{q}}$
denote the theta null point of $A$ with respect to the canonical theta
structure $\Theta_{q}$.
\begin{theorem}
\label{lst}
There exists an $\omega \in R^*$ such that for all
$u,v \in Z_{q}$ one has
\begin{eqnarray}
\label{relateomega}
x_{u+v}x_{u-v} = \omega 
\sum_{z \in Z_2} x^{\sigma}_{ u + z }  x^{\sigma}_{v + z}.
\end{eqnarray}
\end{theorem}
\begin{proof}
We can assume that we have chosen an isomorphism
\begin{eqnarray}
\label{triviaplus1}
Z_{2q,R} \stackrel{\sim}{\rightarrow}
A[2q]^{\mathrm{et}}
\end{eqnarray}
which extends the given trivialization (\ref{trivia}).
We remark that in order to choose the isomorphism
(\ref{triviaplus1}) one may has to extend the base
locally-\'{e}tale. Our assumption is justified by
the fact that the resulting formulas are already defined
over the given ring $R$.
By \cite[Th.2.2]{ca07} there exists a canonical theta
structure $\Theta_{2q}$ for the line bundle $\pol^{2q}$ of
type $Z_{2q}$ which depends only on the
isomorphism (\ref{triviaplus1}).
\newline\indent
Assume that we are given a rigidification of the line bundle $\pol$.
For an abstract finite constant commutative
group $K$ we denote $V(K)= \underline{\mathrm{Hom}}(K_R, \mathcal{O}_R)$.
Let $\pi:A \rightarrow \spec{R}$ denote the structure morphism of the abelian
scheme $A$.
By general theory one can choose theta group equivariant isomorphisms
\[
\mu_{q}: \pi_* \pol^{q} \stackrel{\sim}{\rightarrow}
V( Z_{q}) \quad \mbox{and}
\quad \mu_{2q}: \pi_* \pol^{2q} \stackrel{\sim}{\rightarrow}
V( Z_{2q})
\]
which are uniquely determined up to scalar.
As explained in \cite[$\S$1]{mu66}, together with the chosen
rigidification, these isomorphisms determine finite theta functions $q_{\pol^{q}}$
and $q_{\pol^{2q}}$ in $V( Z_{q})$ and
$V( Z_{2q})$, respectively, which give the coordinates of
the induced theta null points.
It follows from Lemma \ref{comp1} that the theta structures
$\Theta_{q}$ and $\Theta_{2q}$ are $2$-compatible, and hence
satisfy the compatibility assumptions of
Theorem \ref{mult}.
As a consequence, by Corollary \ref{qell2} there exists a $\lambda \in R^*$ such that
\begin{eqnarray}
\label{firstform}
q_{\pol^{q}}(u+v) q_{\pol^{q}}(u-v) = \lambda  \sum_{z \in Z_2} q_{\pol^{2q}}(u+z) q_{\pol^{2q}}(v+z)
\end{eqnarray}
for all $u,v \in Z_q$.
By Theorem \ref{cantwisttheta}, which describes the action of $\sigma$
on the canonical theta null point, by Lemma \ref{comp2} and
by \cite[Lem.3.5]{ckl08}
we deduce from equation
(\ref{firstform}) that
\[
q_{\pol^{q}}(u+v) q_{\pol^{q}}(u-v) = \lambda  \sum_{z \in Z_2} q_{\pol^{q}}^{\sigma}(u+z) q_{\pol^{q}}^{\sigma}(v+z).
\]
This completes the proof of the theorem.
\end{proof}
We remark that one has $x_0 \in R^*$. For the case $q=2$ this follows immediately from the specific shape of the equations (\ref{specialcase}). The case $q>2$ can be reduced to the case $q=2$ by descending
along the relative $\frac{q}{2}$-Frobenius and by applying
\cite[Lem.3.5]{ckl08}.
\newline\indent
Assume now that $k$ is a finite field of
characteristic $2$ and that we have normalized the theta null point
such that $x_0=1$.
The resulting $\omega$ then is is expected to be equal to
the inverse of the determinant of relative Verschiebung
(up to norm one unit).
Evidence is given by the example of
Section \ref{exone}.

\section{Examples of canonical theta null points}
\label{ccft}

In the following we will illustrate Theorem \ref{lst} by some examples.
The following lemma relates the theta null points of
elliptic curves to the
classical $\lambda$- and $j$-invariants.
The lemma will be of use in the following sections.
Let $K$ be a field of characteristic unequal to $2$.
Assume that we are given an elliptic curve $E$ over $K$.
We denote the ample line bundle associated to the zero section
of $E$ by $\pol$. Assume that we are given a theta structure
$\Theta$ of type $\xz / 2 \xz$ for the line bundle $\pol^2$.
We denote the theta null point of $E$ with respect
to the triple $(E,\pol^2,\Theta)$ by $(a_0,a_1)$.
\begin{lemma}
\label{lambda}
The curve $E$ can be given by the equation
\[
y^2=x(x-1)(x-\lambda) \quad \mbox{where} \quad
\lambda = \left( \frac{a_0^2-a_1^2}{a_0^2+a_1^2} \right)^2.
\]
\end{lemma}
\begin{proof}
For simplicity we assume that $(a_0,a_1)=(1,\mu)$ with $\mu \in K$.
The proof in the case $(a_0,a_1)=(\mu,1)$ is analogous.
The morphism $\tau:E \rightarrow \xp^1_{K}$
induced by the theta structure $\Theta$ for $\pol^2$
is surjective of degree $2$.
The group $G(\pol^{\otimes 2}) / \xg_{m,K} \cong A[\pol^2]$
acts on $E$ by translation.
The theta structure $\Theta$ induces a Lagrangian structure
$(\xz / 2 \xz)_{K}
\times \mu_{2,K} \stackrel{\sim}{\rightarrow} A[\pol^2]$.
We claim that
the group elements $(1,1)$ and $(0,-1)$ act on $\xp^1_{K}$ by the matrices
\[
\left( \begin{array}{cc}
0 & 1 \\ 1 & 0 \end{array} \right) \quad \mbox{and}
\quad \left( \begin{array}{cc}
1 & 0 \\ 0 & -1 \end{array} \right).
\]
We set
\[
\delta_z(x)= \left\{ \begin{array}{l@{, \quad}l}
1 & x=z \\
0 & x \not=z
\end{array}
\right.
\]
for $x,z \in \ZZ / 2 \ZZ$. 
Using the definition of the action of the standard theta group
$G(\xz/2\xz)$ on the module of finite theta functions
one computes
\[
(1,x,1) \delta_z = \delta_{z-x} \quad \mbox{and} \quad
(1,0,l) \delta_z = l(z)  \delta_z.
\]
This proves our claim.
By construction the theta null point $(1,\mu)$ induces
a ramification point of the morphism
$\tau$.
The orbit of $(1, \mu)$ under the action of the group $A[\pol^2]$ is given by
$(1,\mu)$, $(1,-\mu)$, $(\mu,1)$ and $(-\mu,1)$.
Clearly, the points in the orbit of $(1,\mu)$
give rise to ramification points of $\tau$.
By Hurwitz's theorem
there are exactly $4$ ramification points of $\tau$.
We map $(1,\mu) \mapsto (0,1)$, $(1,-\mu) \mapsto (1,0)$ and
$(\mu,1) \mapsto (1,1)$ by the linear transformation
\[
\left( \begin{array}{cc} \frac{\mu^2+1}{\mu^2-1} &
- \frac{\mu^2+1}{\mu(\mu^2-1)} \\
1 & \frac{1}{\mu} \end{array} \right).
\]
The latter transformation maps the point $(-\mu,1)$ to $(1,\lambda)$.
This completes the proof of the lemma.
\end{proof}

\subsection{Elliptic curves with $2$-theta structure}
\label{exone}

Let $E$ be an elliptic curve over $\xz_q$ where $\xz_q$ denotes the
Witt vectors with values in a finite field $\xf_q$ with $q=2^d$ elements.
By $\sigma$ we denote the canonical lift of the absolute $2$-Frobenius of $\xf_q$ to
$\xz_q$.
Assume that $E$ has ordinary reduction and that $E$ is
the canonical lift of $E_{\xf_q}$.
Let $\pol$ denote the ample symmetric line bundle associated to
the Weil divisor given by the zero section $0_E$ of $E$.
There exists a unique isomorphism
\[
(\xz / 2 \xz)_{\xz_q} \stackrel{\sim}{\rightarrow} E[2]^{\mathrm{et}}.
\]
By \cite[Th.2.2]{ca07} there exists a
canonical theta structure of type $\xz / 2 \xz$
for the pair $(E,\pol^{\otimes 2})$.
By Theorem \ref{lst} there exists an $\omega \in \xz_q^*$ such that
the coordinates of the theta null point $(x_0,x_1)$ with respect to
to the canonical theta structure satisfy the equations
\begin{eqnarray}
\label{magic}
x_0^2 = \omega  \big( (x_0^{\sigma})^2 +(x_1^{\sigma})^2 \big)
\quad \mbox{and} \quad
x_1^2=2 \omega  x_0^{\sigma} x_1^{\sigma}.
\end{eqnarray}
The equations (\ref{magic}) imply that $x_1 \equiv 0 \bmod 2$.
Hence $x_0 \in \xz_q^*$.
We set $\mu=x_1/x_0$.
Rewriting the equations (\ref{magic}) in terms of $\mu$ we get
\begin{eqnarray}
\label{significant}
\mu^2  \big( (\mu^{\sigma})^2+1 \big)=2 \mu^{\sigma}.
\end{eqnarray}
We set $L=\mathrm{End}_{\xz_q}(E) \otimes \xq$.
\newline\newline\noindent
\bfseries Case $d=1$. \mdseries Equation (\ref{significant}) implies that
\[
0=\mu^3 + \mu - 2 = (\mu-1)  (\mu^2+\mu+2).
\]
A short calculation shows that $L=\xq(\sqrt{-7})$ which has class number $1$.
The polynomial $x^2+x+2$ is reducible over $L$.
We remark that $\mu=2 \omega$.
Note that as expected $\omega$ is the inverse of the invertible
eigenvalue of the $2$-Frobenius endomorphism of $E_{\xf_2}$.
Lemma \ref{lambda} implies that the $j$-invariant of the
elliptic curve $E$, which corresponds to $\mu$, equals
$-15^3$.
\newline\newline\noindent
\bfseries Case $d=2$. \mdseries Equation (\ref{significant}) implies that
\[
0= (\mu^2+\mu+2)  (\mu^4+4 \mu^3
+5 \mu^2 + 2 \mu +4).
\]
Assume that the $j$-invariant of $E_{\xf_4}$ is not equal to $1$.
Then $L=\xq(\sqrt{-15})$ which has class number $2$.
We remark that the minimal polynomial $x^4+4x^3+5x^2+2x+4$
of $\mu$ generates
the Hilbert class field of $L$.
Using the right hand equation (\ref{magic}) and
the relation $\mu^{\sigma^2}=\mu$
one deduces that $\mu^3=8 \omega^2 \omega^{\sigma}$.
We conclude that $\mu=2\zeta\omega$ where $\zeta$ is a norm one unit
in $\xz_q$.
This gives evidence for our conjecture
that $\omega$ is the inverse
of the determinant of the relative Verschiebung morphism (up to norm
one unit). 

\subsection{Elliptic curves with $4$-theta structure}
\label{extwo}

Let $\ZZ_q$, $\QQ_q$, $\sigma$, $E$ and $\pol$ be as in Section \ref{exone}.
Assume that we are given an isomorphism
\begin{eqnarray}
\label{twotrivv}
(\xz / 4 \xz)_R \stackrel{\sim}{\rightarrow} E[4]^{\mathrm{et}}.
\end{eqnarray}
By \cite[Th.2.2]{ca07} there exists a
canonical theta structure of type $(\xz / 4 \xz)$
for the pair $(A,\pol^{\otimes 4})$ depending on the trivialization
(\ref{twotrivv}).
The canonical theta structure induces a closed immersion
$\tau:E \rightarrow \xp^3_{\ZZ_q}$.
Let $(x_0,x_1,x_2,x_3)$
denote the image of the zero section of $E$ under $\tau$.
According to Mumford the image of $\tau_{\QQ_q}$ in $\xp^3_{\QQ_q}$
is the intersection of the quadratic hyper-surfaces
\begin{eqnarray}
\label{mumfordeq}
y_1^2 + y_3^2 = 2 \lambda  y_0  y_2
\quad \mbox{and} \quad y_0^2+y_2^2=2\lambda  y_1  y_3
\end{eqnarray}
where $\lambda = \frac{x_1^2}{x_0x_2}$.
By symmetry the theta null point lies in the plane $y_1=y_3$.
For more details see \cite[$\S$5]{mu66}.
We can assume that $x_0=1$.
By Theorem \ref{lst} there exists a unique $\omega \in R^*$ such that
\begin{eqnarray}
\label{eeins}
& & 1 = \omega  \big( 1 + (x_2^{\sigma})^2 \big) \\
\label{ezwei}
& & x_1^2 = \omega  x_1^{\sigma} (1+x_2^{\sigma} ) \\
\label{edrei}
& & x_2^2 = 2 \omega  x_2^{\sigma} \\
\label{evier}
& & x_2= 2\omega (x_1^{\sigma})^2.
\end{eqnarray}
Now assume that $d=1$.
By equation (\ref{edrei}) we have $x_2=2 \omega$.
It follows by the equation (\ref{ezwei}) that $x_1=x_3=\omega(1+2\omega)$.
As a consequence we have
\[
\lambda= \frac{\omega}{2}  (1+ 2 \omega)^2
\]
where $\lambda$ is as above.
Finally we conclude by equation (\ref{eeins}) that
\[
4\omega^3+\omega-1= \big( 2 \omega-1 \big) \big (2\omega^2+\omega
+ 1 \big)= 0.
\]
The latter factor makes clear that $\omega$ is a reciprocal 
eigenvalue of Frobenius in an order of discriminant $-7$.

\section*{Acknowledgments}

The results of this article were proved with some help of
B.~Edixhoven.
I am grateful to D.~Kohel for some very valuable comments regarding
the examples of Section \ref{ccft}.

\bibliographystyle{plain}

\vspace{1cm}
\noindent
\begin{tabular}{l}
Robert Carls \\
\texttt{robert.carls@uni-ulm.de} \\ 
Institute of Pure Mathematics \\
University of Ulm \\
D-89069 Ulm, Germany
\end{tabular}

\end{document}